\documentclass[a4paper, 11pt, twoside]{article}
\usepackage{mfpaperstuff2}
\begin{document}
\newcommand\see{\diamondsuit}
\newcommand\se[2]{\see_{#1}(#2)}
\newcommand\sof[2]{#1^{(#2)}}
\newcommand\smx{sometimes maximal\xspace}
\newcommand\smn{sometimes minimal\xspace}
\newcommand\mcd{multipartition datum\xspace}
\newcommand\mcds{multipartition data\xspace}
\newcommand\hk[2]{\operatorname{hk}_{#1}(#2)}
\renewcommand\gcd{greatest common divisor\xspace}
\newcommand\gcds{greatest common divisors\xspace}
\newcommand\ngcd{\operatorname{gcd}}
\newcommand\ltup[2]{\left(\left.#1\ \right|\ \smash{#2}\right)}
\newcommand\rtup[2]{\left(\smash{#1}\ \left|\ #2\right.\right)}
\newcommand\weyl{\tilde{\mathfrak S}_}
\newcommand\be[1]{\mathcal{B}^{#1}}
\newcommand\ber[2]{\be{#1}_{#2}}
\newcommand\berm[2]{\be{#1}_{-#2}}
\newcommand\cores[1]{\calc_{#1}}
\newcommand\zsz{\bbz/s\bbz}
\newcommand\zgz{\bbz/g\bbz}
\newcommand\bc[2]{\mc{#1}{0,#2}}
\newcommand\mc[2]{\mathcal{C}_{(#1\,|\,#2)}}
\newcommand\smc[1]{\mathcal{C}_{#1}}
\newcommand\cpt[2]{#1^{(#2)}}
\Yvcentermath1

\title{Simultaneous core multipartitions}

\msc{05A17, 05E10}
\toptitle

\begin{abstract}
We initiate the study of \emph{simultaneous core multipartitions}, generalising simultaneous core partitions, which have been studied extensively in the recent literature. Given a \emph{multipartition datum} $(s\,|\,c)$, which consists of a non-negative integer $s$ and an $l$-tuple $c$ of integers, we introduce the notion of an $(s\,|\,c)$-core multipartition. Given an arbitrary set of \mcds, we give necessary and sufficient conditions for the corresponding set of simultaneous core multipartitions to be finite. We then study the special case of simultaneous core bipartitions, giving exact enumerative results in some special subcases.
\end{abstract}

\section{Introduction}

The study of integer partitions has a long history, with applications in a variety of areas. In the last few years there has been considerable interest in \emph{core} partitions, i.e.\ partitions with no hook lengths divisible by a certain prescribed integer. Of particular interest are \emph{simultaneous core partitions}, i.e.\ partitions which are both $s$- and $t$-cores, for given (typically coprime) integers $s,t$. Various enumerative results have been proved for these ``$(s,t)$-cores''; foremost among these are Anderson's theorem \cite{and} giving the number of $(s,t)$-cores, and Armstrong's Conjecture (stated in \cite{ahj}, and proved by Johnson in \cite{j}) giving the average size of an $(s,t)$-core.

In this paper we introduce the subject of core multipartitions. For a fixed $l\in\bbn$, an $l$-multipartition is just an $l$-tuple of partitions. We generalise the notion of core partition to multipartitions by using a characterisation of core partitions in terms of residues of nodes which goes back to a result of Littlewood \cite{litt}. Our definition of core multipartitions has representation-theoretic significance in terms of modules for cyclotomic Hecke algebras. The extension to multipartitions requires not just an integer $s$ but also an $l$-tuple $c\in\bbz^l$; so we actually introduce the notion of an \emph{$(s\,|\,c)$-core multipartition} (or simply an $(s\,|\,c)$-core). We can then consider the question of simultaneous core multipartitions, i.e.\ multipartitions which are $(s\,|\,c)$-cores for all pairs $(s\,|\,c)$ in a given set $\calt$. Our main result is a determination of exactly when there are only finitely such multipartitions.

Having established this finiteness result, we consider enumerative results, restricting to the first non-trivial case (where $l=2$ and $\cardx\calt=2$) and enumerating simultaneous core bipartitions in two special subcases.

\section{Definitions and basic results}

%In this section we set out some definitions and notation.

\subsection{Standard notation}

If $X\subseteq\bbz^l$, $n\in\bbz^l$ and $s\in\bbn$, then we write $n+sX=\lset{n+sx}{x\in X}$. We define $\zsz$ to be the set $\lset{n+s\bbz}{n\in\bbz}$. (We do not employ the popular abuse of notation in which $\zsz=\{0\dots,s-1\}$.) A \emph{$\zsz$-tuple} of integers just means a function $u$ from $\zsz$ to $\bbz$, which we write in the form $\ltup{u_i}{i\in\zsz}$.

\subsection{Partitions}

A \emph{partition} is a weakly decreasing sequence $\la=(\la_1,\la_2,\dots)$ of non-negative integers with finite sum. When writing partitions, we typically group together equal parts with a superscript and omit the trailing zeroes, and we write the partition $(0,0,\dots)$ as $\varnothing$. We let $\calp$ denote the set of all partitions.

The \emph{size} of a partition $\la$ is the sum $\cardx\la=\sum_{a\gs1}\la_a$. The \emph{Young diagram} of $\la$ is the set
\[
[\la]=\lset{(a,b)\in\bbn^2}{b\ls\la_a}
\]
whose elements we call the \emph{nodes} of $\la$. We draw $[\la]$ as an array of boxes in the plane using the English convention, in which the Young diagram of $(6,4,2,1^2)$ is drawn as follows.
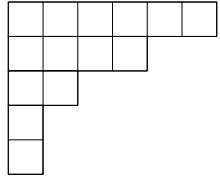
\[
\yng(6,4,2,1^2)
\]
A node of $\la$ is \emph{removable} if it can be removed to leave a Young diagram (i.e.\ if it has the form $(a,\la_a)$, with $\la_a>\la_{a+1}$), while a pair $(a,b)\notin[\la]$ is an \emph{addable node} of $[\la]$ if it can be added to $[\la]$ to yield a Young diagram. If $(a,b)$ is a node of $\la$, the \emph{$(a,b)$-hook} of $\la$ is the set of nodes of $\la$ directly to the right of $(a,b)$ or directly below $(a,b)$, including $(a,b)$ itself. The \emph{$(a,b)$-hook length} is the number of nodes in this hook. If the $(a,b)$-hook has length $s$, we call it an \emph{$s$-hook}. $\la$ is an \emph{$s$-core partition} (or simply an \emph{$s$-core}) if it has no $s$-hooks.

For example, the shaded nodes in the diagram below comprise a $5$-hook of $(6,4,2,1^2)$, so this partition is not a $5$-core. On the other hand, one can easily check that this partition has no $3$-hooks, so is a $3$-core.
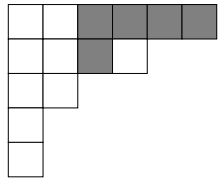
\[
\gyoung(;;!<\Yfillcolour{gray}>;;;;,!<\Yfillcolour{white}>;;!<\Yfillcolour{gray}>;!<\Yfillcolour{white}>;,;;,;,;)
\]
We write $\cores s$ for the set of all $s$-core partitions. These partitions can also be characterised in terms of residues of nodes. Given a node $(a,b)\in[\la]$, define its \emph{$s$-residue} to be $b-a+s\bbz$. For example, the $3$-residues of the nodes of $(6,4,2,1^2)$ are illustrated in the following diagram (in which we label a node of residue $i+3\bbz$ with $i$, for $i=0,1,2$).
\[
\yngres(3,6,4,2,1^2)
\]
The \emph{$s$-content} of a partition $\la$ is defined to be the multiset of $s$-residues of the nodes of $\la$. From the diagram above, we see that the $3$-content of $(6,4,2,1^2)$ is $\{(0+3\bbz)^4,(1+3\bbz)^4,(2+3\bbz)^6\}$ (where we adopt what we hope is an obvious notation for a multiset of elements of $\bbz/s\bbz$).

The $s$-content of a partition is significant because of the following result.

\begin{thmciting}\textup{\textbf{\cite{litt}.}}\label{rescore}
Suppose $\la\in\calp$ and $s>1$. Then $\la$ is an $s$-core \iff there is no other partition with the same $s$-content as $\la$.
\end{thmciting}

Now suppose $s_1,\dots,s_r\in\bbn$. An \emph{$(s_1,\dots,s_r)$-core} means a partition which is an $s_i$-core for each $i$. It is fairly easy to show that the set of $(s_1,\dots,s_r$)-cores is finite \iff $s_1,\dots,s_r$ are coprime; this appears to have been written down for the first time by Xiong \cite[Theorem 1.1]{xiong}. Our main aim in this paper is to prove an analogue of this statement for multipartitions, which we introduce next.

\subsection{Multipartitions}\label{multipintrosec}

%Now we consider multipartitions. 
Fix $l\in\bbn$. An \emph{$l$-multipartition} is an $l$-tuple $\la=(\la^{(1)},\dots,\la^{(l)})$ of partitions, which we call the \emph{components} of $\la$. We write $\calp^l$ for the set of all $l$-multipartitions, and we write $\varnothing^l$ for the multipartition $(\varnothing,\dots,\varnothing)$.

The size of an $l$-multipartition $\la$ is the sum of the sizes of its components. The Young diagram of $\la$ is the set
\[
[\la]=\lset{(a,b,k)\in\bbn^2\times\{1,\dots,l\}}{b\ls\la^{(k)}_a},
\]
whose elements we call the nodes of $\la$. We draw the Young diagram of $\la$ by drawing the Young diagrams of $\la^{(1)},\dots,\la^{(l)}$ in order from left to right. We define addable and removable nodes of multipartitions analogously to those for partitions.

Now take an $l$-tuple $c=(c_1,\dots,c_l)\in\bbz^l$, and define the \emph{$(s\,|\,c)$-residue} of a node $(a,b,k)$ to be $b-a+c_k+s\bbz$. We refer to a node of $(s\,|\,c)$-residue $i\in\zsz$ as an \emph{$i$-node}. Define the \emph{$(s\,|\,c)$-content} of $\la$ to be the multiset of $(s\,|\,c)$-residues of the nodes of $\la$.

For example, suppose $l=3$, $s=4$ and $c=(0,2,1)$. For $\la=\left((2),(4,1^2),(1^2)\right)$, the residues are indicated by the following diagram.
\Yvcentermath0
\[
\yngres(4,2)\qquad\Ycorner2\yngres(4,4,1^2)\qquad\Ycorner1\yngres(4,1^2)
\]
We see that the $(4\,|\,(0,2,1))$-content of $\la$ is $\left\{(0+4\bbz)^4,(1+4\bbz)^4,(2+4\bbz)^1,(3+4\bbz)^1\right\}$.

Now, inspired by \cref{rescore}, we make the following definition: say that $\la\in\calp^l$ is an \emph{$(s\,|\,c)$-core multipartition} (or simply an \emph{$(s\,|\,c)$-core}) if there is no other $l$-multipartition $\mu$ with the same $(s\,|\,c)$-content. We write $\mc sc$ for the set of all $(s\,|\,c)$-cores.

In the rest of the paper we will refer to the pair $(s\,|\,c)$ as an \emph{$l$-\mcd}. If we wish to refer to the individual integers $c_1,\dots,c_l$, we may write $(s\,|\,c_1,\dots,c_l)$ instead of $(s\,|\,c)$ in any of the notation introduced above. For example, we can easily check that the multipartition $\la$ above lies in $\mc4{0,2,1}$.

\begin{rmksenum}
%\item
%The definition of the $(s\,|\,c)$-residue of a node (and hence of an $(s\,|\,c)$-core) depends only on the residues of $c_1,\dots,c_l$ modulo $s$; so we could just as easily take $c\in(\bbz/s\bbz)^l$ rather than $c\in\bbz^l$ in the definition. It will be slightly more convenient to take $c_1,\dots,c_l$ to be integers. Note also that the set $\mc sc$ is unchanged if $c_1,\dots,c_l$ are all translated by a fixed integer (since the effect is just to shift the residues of all nodes by this same integer).
\item
In the study of $s$-core partitions, $s$ is typically assumed to be greater than~$1$. However, one can meaningfully consider the cases $s=0$ and $s=1$. First take $s=1$; according to the definition using hooks, the only $1$-core partition is $\varnothing$. However, the partition $(1)$ is the unique partition with its $1$-content (which is why we need to assume $s\neq1$ in \cref{rescore}). Nevertheless, most of the theory of core partitions applies (in a trivial way) with $s=1$, if we take $\cores1=\{\varnothing\}$. Similarly for multipartitions, we take $\mc1c=\{\varnothing^l\}$ (which is consistent with the definition of core multipartitions given above provided $l\gs2$), and the results we prove below will apply in this case.

Now consider the case $s=0$. In this case we should regard the residue of a node $(a,b)$ as the integer $b-a$, and correspondingly define the $0$-content to be a multiset of integers. It is then not hard to prove that any partition is determined by its $0$-content, so every partition is a $0$-core.% In fact, the case $s=0$ can be regarded as the limit of taking $s$ very large (in fact, this case is often described as $s=\infty$ rather than $s=0$); since any partition is an $s$-core for all sufficiently large $s$, it is reasonable that every partition should be a $0$-core.

The situation with multipartitions is less straightforward when $s=0$.  Given $c\in\bbz^l$ we define the $(0\,|\,c)$-residue of the node $(a,b,k)$ to be the integer $b-a+c_k$, so that the $(0\,|\,c)$-content is again a multiset of integers. But now not every multipartition is a $(0\,|\,c)$-core. In fact this is easily seen: if $c_j=c_k$ for some $j\neq k$, then a multipartition $\la$ has the same $(0\,|\,c)$-content as the multipartition obtained by switching the components $\cpt\la j$ and $\cpt\la k$, so cannot be a $(0\,|\,c)$-core if these components are unequal. So the study of $(0\,|\,c)$-core multipartitions is certainly non-trivial, and we will include the case $s=0$ in our considerations in this paper. Given two integers $a,b$, the condition $a\equiv b\ppmod s$ should be read as $a=b$ in the case $s=0$.

The case $s=0$ can be regarded as the limiting case as $s$ gets very large (in fact, this situation is often described as $s=\infty$ rather than $s=0$): given a multipartition $\la$ and $c\in\bbz^l$, it is easily seen that we have $\la\in\mc0c$ \iff $\la\in\mc sc$ for all sufficiently large $s$.
\item
Our definition of core multipartitions is not completely arbitrary, but has representation-theoretic significance. Associated to a pair $(s\,|\,c)$ as above and a positive integer $n$ is an \emph{Ariki--Koike algebra} (a Hecke algebra of the complex reflection group of type $G(l,1,n)$). This algebra has an important family of modules (the \emph{Specht modules}) labelled by $l$-multi\-partitions of size $n$. Lyle \& Mathas \cite{lyma} showed that two multipartitions having the same $(s\,|\,c)$-content is equivalent to the corresponding Specht modules lying in the same block of the Ariki--Koike algebra, and the author \cite{mfwt} showed that a multipartition being an $(s\,|\,c)$-core is equivalent to the corresponding Specht module being contained in a simple block. This is analogous to the significance of $s$-cores in the $s$-modular representation theory of the symmetric group (or more generally the Iwahori--Hecke algebra of type $A$). In fact, this analogy goes further: in \cite{mfwt} the author defines a ``weight'' function on multipartitions (depending on $s,c$) which is an analogue of the $s$-weight of a partition $\la$  (i.e.\ the number of rim $s$-hooks that need to be removed to reach the $s$-core of $\la$). $(s\,|\,c)$-cores are then simply multipartitions of weight $0$. We will use some of the results from \cite{mfwt} below.
\end{rmksenum}

\subsection{Basic results}\label{basicmcoresec}

In this section we will give some basic results on core multipartitions; in particular, we will give a simple condition in terms of beta-numbers for a multipartition to be an $(s\,|\,c)$-core.

We start with two very simple results.

\begin{lemma}\label{shiftmods}
Suppose $(s\,|\,c)$ and $(s\,|\,d)$ are two $l$-\mcds, with $d_i-c_i\equiv d_j-c_j\ppmod s$ for all $i,j$. Then $\mc sc=\mc sd$.
\end{lemma}

\begin{pf}
First note that the set $\mc sc$ is unchanged if we add a fixed integer $a$ to each $c_i$, since the effect is just to shift the residues of all nodes by $a$. Doing this with $a=d_1-c_1$, we may assume that $c_i\equiv d_i\ppmod s$ for each $i$. But then the $(s\,|\,c)$-residue of each node is the same as the $(s\,|\,d)$-residue, so that the $(s\,|\,c)$-content of any multipartition is the same as the $(s\,|\,d)$-content, and the result follows.
\end{pf}

\begin{lemma}\label{multicore}
Suppose $(s\,|\,c)$ is a \mcd, and $\la\in\mc sc$. Then each component of $\la$ is an $s$-core.
\end{lemma}

\begin{pf}
We prove the contrapositive. Suppose $\cpt\la k$ is not an $s$-core. Then by \cref{rescore} there is another partition $\mu$ with the same $s$-content as $\cpt\la k$. The multipartition obtained from $\la$ by replacing $\cpt\la k$ with $\mu$ then has the same $(s\,|\,c)$-content as $\la$, so $\la$ is not an $(s\,|\,c)$-core.
\end{pf}

Our remaining background results are mostly taken from \cite{mfwt}, but we need to explain how to translate the results into our notation.

The combinatorics in \cite{mfwt} are based on data consisting of a field $\bbf$ and non-zero elements $q,Q_1,\dots,Q_l$ of $\bbf$. (In fact the integer $r$ is used instead of $l$ in \cite{mfwt}, but this makes no practical difference.) The residue of a node $(a,b,k)$ is defined in \cite{mfwt} to be the element $q^{b-a}Q_k$ of $\bbf$. For the purposes of the present paper, it suffices to assume that each $Q_i$ is a power of $q$, say $Q_i=q^{c_i}$, for $c_i\in\bbz$. If we let $s$ denote the multiplicative order of $q\in\bbf$ (this order is called $e$ in \cite{mfwt}), then two nodes have the same residue (in the sense of \cite{mfwt}) \iff they have the same $(s\,|\,c_1,\dots,c_l)$-residue. (The $(s\,|\,c_1,\dots,c_l)$-residue is essentially the base $q$ logarithm of the residue in \cite{mfwt}). Given a multipartition $\la$ and $f\in\bbf$, \cite{mfwt} defines $c_f(\la)$ to be the number of nodes of $\la$ of residue $f$, and defines the \emph{weight} of $\la$ to be
\[
w(\la)=\sum_{i=1}^lc_{Q_i}(\la)-\frac12\sum_{f\in\bbf}\left(c_f(\la)-c_{qf}(\la)\right)^2.
\]
Multipartitions $\la,\mu$ are defined to lie in the same \emph{combinatorial block} \iff $c_f(\la)=c_f(\mu)$ for every $f$. Clearly, this is equivalent to $\la$ and $\mu$ having the same $(s\,|\,c_1,\dots,c_l)$-content. Thus a multipartition $\la$ is an $(s\,|\,c_1,\dots,c_l)$-core \iff it lies in a combinatorial block by itself. \cite[Theorem 4.1]{mfwt} then says that this happens \iff $w(\la)=0$.

The results in \cite[Section~3]{mfwt} yield a simple algorithm for computing the weight of a multipartition, and in particular for determining whether a multipartition has weight $0$. The first result that we cite shows that in order to check whether a multipartition is an $(s\,|\,c)$-core we can reduce to the case $l=2$.

\begin{propn}\label{lev2}
Suppose $(s\,|\,c)$ is an $l$-\mcd, and that $\cpt\la k\in\cores s$ for every $k$. Then $\la$ is an $(s\,|\,c)$-core \iff $(\cpt\la j,\cpt\la k)$ is an $(s\,|\,c_j,c_k)$-core bipartition for all $1\ls j<k\ls l$.
\end{propn}

\begin{pf}
Define the weight function $w$ as above using the data $q,Q_1,\dots,Q_l$, where $q$ is a primitive $s$th root of unity (or a non-zero non-root of unity, if $s=0$) in a field $\bbf$, and $Q_i=q^{c_i}$ for each $i$. \cite[Proposition 3.5]{mfwt} says (given the assumption that each $\cpt\la k$ lies in $\cores s$) that $w(\la)$ is the sum of the values $w((\la^{(j)},\la^{(k)}))$ over all pairs $j<k$, where $w((\la^{(j)},\la^{(k)}))$ is defined using the data $q,Q_j,Q_k$. Since weight is non-negative by \cite[Corollary 3.9]{mfwt}, this means that $w(\la)=0$ (i.e.\ $\la$ is an $(s\,|\,c)$-core) \iff $w((\la^{(j)},\la^{(k)}))=0$ for every $j,k$ (i.e.\ each $(\la^{(j)},\la^{(k)})$ is an $(s\,|\,c_j,c_k)$-core).
\end{pf}

To go further, we recall the definition of beta-numbers, which goes back to Nakayama \cite{nak}. Define the \emph{beta-set} of a partition $\la$ to be the set
\[
\be\la=\lset{\la_a-a}{a\in\bbn}.
\]
For any $c\in\bbz$, we write $\ber\la c$ for the set $\be\la+c$, which we refer to as the \emph{$c$-shifted beta-set} of $\la$.

The following result is due to Robinson \cite[(2.8)]{robin4}.

\begin{propn}\label{corebeta}
Suppose $s\gs0$ and $\la\in\calp$. Then the number of $s$-hooks of $\la$ equals the number of $b\in\be\la$ such that $b-s\notin\be\la$. In particular, $\la$ is an $s$-core \iff $\be\la\supseteq\berm\la s$.
\end{propn}

This result is key in the study of core partitions; it yields James's abacus model \cite[Section~2.7]{jk} for partitions, which in turn leads to a geometric interpretation for the set of $s$-cores.

We make an observation about beta-sets which will be useful later. Suppose $\la,\mu\in\calp$ and $c,d\in\bbz$. Note that $\ber\la c$ is a set of integers which is bounded above and whose complement in $\bbz$ is bounded below. Moreover, the number of non-negative integers in $\ber\la c$ minus the number of negative integers not in $\ber\la c$ equals $c$. As a consequence, we have the following result.

\begin{lemma}\label{betaobs}
Suppose $c,d\in\bbz$ and $\la,\mu\in\calp$. Then $\card{\ber\la c\setminus\ber\mu d}-\card{\ber\mu d\setminus\ber\la c}=c-d$. In particular, if $\ber\la c\supseteq\ber\mu d$, then $c\gs d$.
\end{lemma}

We now explain how core multipartitions can be characterised in terms of the beta-sets of their components. In view of \cref{lev2} we restrict to the case $l=2$.% The next two results are special cases of the results in \cite[\S3.4]{mfwt}.

\begin{propn}\label{sbetacore}
Suppose $c,d\in\bbz$ and $(\la,\mu)\in\calp^2$.
\begin{enumerate}[ref=\cref{sbetacore}(\arabic*)]
\item\label{s0betacore}
\begin{enumerate}
\item
If $c\ls d$, then $(\la,\mu)\in\mc0{c,d}$ \iff $\ber\la c\subseteq\ber\mu d$.
\item
If $c\gs d$, then $(\la,\mu)\in\mc0{c,d}$ \iff $\ber\la c\supseteq\ber\mu d$.
\end{enumerate}
\item\label{sposbetacore}
Suppose $s\in\bbn$, and let $e$ be the residue of $c-d$ modulo $s$. Suppose $(\la,\mu)\in\calp^2$. Then $(\la,\mu)\in\mc s{c,d}$ \iff
\[
\ber\la e\supseteq\be\mu\supseteq\ber\la{e-s}.
\]
\end{enumerate}
\end{propn}

\begin{pf}
Following \cite{mfwt} we define integers $\gamma_i$ for $i\in\bbz$ as follows.
\begin{itemize}
\item
If $s=0$, then we set
\[
\gamma_i=
\begin{cases}
1&\text{if }\ber\la c\ni i\notin\ber\mu d\\
-1&\text{if }\ber\la c\not\ni i\in\ber\mu d\\
0&\text{otherwise}.
\end{cases}
\]
\item
If $s>0$, then we define $\gamma_i$ to be the largest element of $\ber\la c\cap(i+s\bbz)$ minus the largest element of $\ber\mu d\cap(i+s\bbz)$.
\end{itemize}

In either case, \cite[Lemma 3.7(3) \& Proposition 3.8]{mfwt} show that $(\la,\mu)$ has weight $0$ (i.e.\ is an $(s\,|\,c,d)$-core) \iff $\gamma_i-\gamma_j\ls1$ for every $i,j$.

If the case $s=0$, the above condition says that $(\la,\mu)$ is an $(s\,|\,c,d)$-core \iff the $\gamma_i$ are all non-negative or all non-positive. But note that
\[
\sum_{i\in\bbz}\gamma_i=\card{\ber\la c\setminus\ber\mu d}-\card{\ber\mu d\setminus\ber\la c}=c-d
\]
by \cref{betaobs}. So if $c\gs d$ then $(\la,\mu)$ is an $(s\,|\,c,d)$-core \iff each $\gamma_i$ is non-negative, which is the same as saying $\ber\la c\supseteq\ber\mu d$. A similar statement applies when $c\ls d$.

If instead $s>0$, then
\[
\sum_{i=0}^{s-1}\gamma_i=\card{\ber\la c\setminus\ber\mu d}-\card{\ber\mu d\setminus\ber\la c}=c-d,
\]
which means we have $\gamma_i-\gamma_j\ls1$ for all $i,j$ \iff $\gamma_i\in\{(c-d-e)/s,1+(c-d-e)/s\}$ for all $i$. The condition that $\gamma_i\gs(c-d-e)/s$ for all $i$ is equivalent to the condition $\ber\la e\supseteq\be\mu$, while the condition that $\gamma_i\ls1+(c-d-e)/s$ for all $i$ is equivalent to $\be\mu\supseteq\ber\la{e-s}$.
\end{pf}

\subsection{Action of the affine symmetric group}

One of the most interesting and useful features of the set of $s$-cores is that it admits a natural action of the affine symmetric group. In this section, we show how this generalises to core multipartitions. This provides a natural proof of the fact that (provided $s\neq1$) the set $\mc sc$ is infinite.

Take $s\gs2$. Recall that the \emph{affine symmetric group} $\weyl s$ is the group of all permutations $g$ of $\bbz$ with the properties that
\begin{itemize}
\item
$g(n+s)=g(n)+s$ for all $n\in\bbz$, and
\item
$g(0)+g(1)+\dots+g(s-1)=0+1+\dots+s-1$.
\end{itemize}
Then $\weyl s$ is a Coxeter group, with generating set $\lset{s_i}{i\in\zsz}$ defined by
\[
s_i(n)=
\begin{cases}
n+1&(n\in i-1)\\
n-1&(n\in i)\\
n&(n\notin i-1,i)
\end{cases}
\]
for $i\in\zsz$ and $n\in\bbz$. The subgroup $\weyl s^0$ generated by $\lset{s_i}{i\neq s\bbz}$ is naturally isomorphic to the symmetric group $\sss s$.

Now suppose $\la$ is a partition and $i\in\zsz$. Define $s_i(\la)$ to be the partition obtained by simultaneously adding all the addable $i$-nodes to $\la$ and removing all the removable $i$-nodes. This defines an action of $\weyl s$ on the set of all partitions. Moreover, the set $\cores s$ is an orbit for this action, so we have a transitive action of $\weyl s$ on $\cores s$. This action was first studied by Lascoux \cite{las}, who showed that the stabiliser of the empty partition is the subgroup $\weyl s^0$ defined above, so that $s$-cores are naturally in bijection with left cosets of $\weyl s^0$ in $\weyl s$.

Next suppose we shift all residues of nodes in $\bbn^2$ by some fixed amount $c$; that is, we redefine the residue of a node $(a,b)$ to be $b-a+c+s\bbz$. Then we can define another action of $\weyl s$ on $\cores s$ in exactly the same way as defined above; this just amounts to twisting Lascoux's action by the automorphism of $\weyl s$ defined by $s_i\mapsto s_{i+c}$ for all $i$. We call this the \emph{$c$-shifted action} of $\weyl s$ on $\cores s$. The stabiliser of $\varnothing$ under the $c$-shifted action is the parabolic subgroup $\weyl s^c$ of $\weyl s$ generated by $\lset{s_i}{i\neq c+s\bbz}$.

Now we consider multipartitions. Suppose we have an $l$-\mcd $(s\,|\,c)$; for the moment we will continue to assume that $s\gs2$ (we will comment below on the case $s=0$). We can define an action of $\weyl s$ on $\calp^l$ analogously to the action on $\calp$ above: if $\la\in\calp^l$, then $s_i(\la)$ is the multipartition obtained by adding all addable $i$-nodes and removing all removable $i$-nodes. Then we have the following.

\begin{figure}[p]
\Yboxdim{7pt}
\newcommand\edg[1]{node[midway,fill=white,inner sep=1pt]{\small$#1$}}
\newcommand\egd[3]{\draw[<->](#1)--(#2)node[midway,fill=white,inner sep=1pt]{\small$#3$};}
\[
\begin{tikzpicture}[scale=2.6]
\draw(0,0)node(1){$\varnothing\ \varnothing$}edge[{}-{},loop left]node{\small$2$}()
++(60:1)node(2){$\varnothing\ \yng(1)$}
++(120:1)node(3){$\varnothing\ \yng(2)$}edge[{}-{},loop left]node{\small$1$}()
++(60:1)node(4){$\yng(1)\ \yng(3,1)$}
++(1,0)node(5){$\yng(1^2)\ \yng(3,1^2)$}
++(300:1)node(6){$\yng(1^2)\ \yng(2,1^2)$}
++(1,0)node(7){$\yng(2,1^2)\ \yng(2^2,1^2)$}
++(300:1)node(8){$\yng(3,1^2)\ \yng(2^2,1^2)$}
++(240:1)node(9){$\yng(3,1)\ \yng(2,1^2)$}
++(300:1)node(10){$\yng(4,2)\ \yng(3,1^2)$}
++(240:1)node(11){$\yng(4,2)\ \yng(3,1)$}
++(-1,0)node(12){$\yng(3,1)\ \yng(2)$}
++(240:1)node(13){$\yng(3,1^2)\ \yng(2)$}
++(-1,0)node(14){$\yng(2,1^2)\ \yng(1)$}
++(120:1)node(15){$\yng(1^2)\ \varnothing$}edge[{}-{},loop left]node{\small$0$}()
++(60:1)node(16){$\yng(1)\ \varnothing$}
++(1,0)node(17){$\yng(2)\ \yng(1)$}
++(60:1)node(18){$\yng(2)\ \yng(1^2)$}
++(120:1)node(19){$\yng(1)\ \yng(1^2)$};
\egd232
%\draw[<->](2)--(3)\edg2;
\egd340
\egd452
\egd560
\egd671
\egd782
\egd891
\egd9{10}0
\egd{10}{11}2
\egd{11}{12}0
\egd{12}{13}1
\egd{13}{14}2
\egd{14}{15}1
\egd{15}{16}2
\egd{16}{17}1
\egd{17}{18}0
\egd{18}{19}1
\egd{19}20
\egd211
\egd1{16}0
\egd6{19}2
\egd{18}92
\egd{17}{12}2
\draw[dashed](4)--++(120:.5);
\draw[dashed](5)--++(60:.5);
\draw[dashed](7)--++(60:.5);
\draw[dashed](8)--++(.5,0);
\draw[dashed](10)--++(.5,0);
\draw[dashed](11)--++(300:.5);
\draw[dashed](13)--++(300:.5);
\draw[dashed](14)--++(240:.5);
\end{tikzpicture}
\]
\caption
{The action of the affine symmetric group of degree $3$ on $\mc3{0,1}$}\label{actfig}
\end{figure}
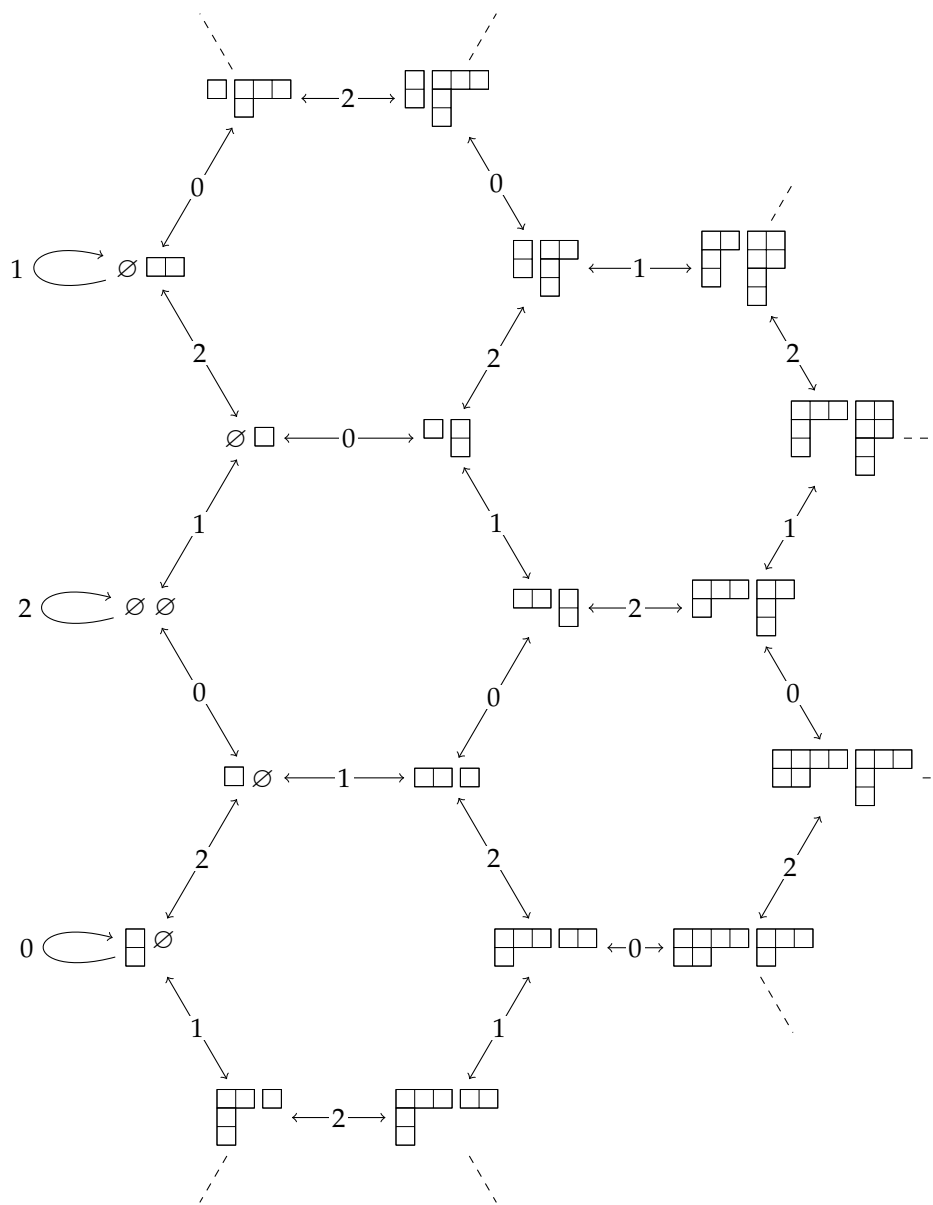

\begin{propn}\label{mcorbit}
Suppose $s\gs2$. Under the action of $\weyl s$ on $\calp^l$ described above, $\mc sc$ is an orbit.
\end{propn}

\begin{pf}
First we show that if $\la\in\mc sc$ and $i\in\zsz$ then $s_i(\la)\in\mc sc$. Note that $\la$ cannot have both addable and removable $i$-nodes, because if it did, then we could remove a removable $i$-node and add an addable $i$-node to obtain another multipartition with the same $(s\,|\,c)$-content, contradicting the assumption that $\la\in\mc sc$. So we assume that $\la$ has no addable $i$-nodes (the other case is similar). Then $s_i(\la)$ is obtained just by removing all the removable $i$-nodes from $\la$. Now by \cite[Lemma 3.6]{mfwt} $\la$ and $s_i(\la)$ have the same weight (note that the integers $u$ and $\delta_i(\la)$ appearing in that lemma are both equal to the number of removable $i$-nodes of $\la$ in our situation, so the term on the right-hand side is zero) and hence $s_i(\la)$ is also an $(s\,|\,c)$-core.

So $\mc sc$ is a union of orbits. To show that $\mc sc$ is a single orbit, we show that if $\la\in\mc sc$ with $\la\neq\varnothing^l$, then there is a strictly smaller multipartition in the same orbit; applying this repeatedly, we find that $\varnothing^l$ lies in the same orbit as $\la$.

The assumption that $\la\neq\varnothing^l$ mean that $\la$ has at least one removable node, of residue $i$, say. As observed at the start of the proof, $\la$ cannot have any addable $i$-nodes, so $s_i(\la)$ is obtained from $\la$ by removing $i$-nodes only. So $s_i(\la)$ is strictly smaller than $\la$, as required.
\end{pf}

Of course, \cref{mcorbit} can be used as an alternative definition of $\mc sc$ in the case $s\neq1$: we can define $\mc sc$ to be the orbit containing $\varnothing^l$ under the action of $\weyl s$ on $\calp^l$.

Part of the action of $\weyl s$ on $\mc sc$ is illustrated in \cref{actfig} in the case $s=3$ and $c=(0,1)$. In this diagram an arrow labelled $i$ indicates the action of $s_{i+3\bbz}$.

In order to understand the action of $\weyl s$ on $\mc sc$ in general, we find the stabiliser of $\varnothing^l$. This is easy to work out, given the discussion above of the shifted actions of $\weyl s$ on $\cores s$. It is clear from the definitions that $g\in\weyl s$ fixes $\varnothing^l$ \iff it fixes $\varnothing$ under the $c_k$-shifted action of $\weyl s$ on $\cores s$, for $k=1,\dots,l$. Hence the stabiliser of $\varnothing^l$ is the intersection $\weyl s^{c_1}\cap\dots\cap\weyl s^{c_l}$. It is a standard fact in the theory of Coxeter groups that the intersection of a family of parabolic subgroups is the parabolic subgroup generated by the intersection of the generating sets of these subgroups. So the stabiliser of $\varnothing^l$ is the subgroup $\lspan{s_i}{i\notin\{c_1+s\bbz,\dots,c_l+s\bbz\}}$. Hence the set $\mc sc$ is in bijection with the set of left cosets of this subgroup.

We now consider the case $s=0$. Here the discussion above applies, except that the finitely-generated Coxeter group $\weyl s$ is replaced with the \emph{finitary symmetric group}, i.e.\ the group $\sss\infty$ of all finitely-supported permutations of $\bbz$. This is also a Coxeter group, with infinite generating set $\lset{s_i}{i\in\bbz}$, where $s_i$ is the transposition $(i-1,i)$. The stabiliser of $\varnothing^l$ under the action of $\sss\infty$ on $\mc0c$ is $\lspan{s_i}{i\notin\{c_1,\dots,c_l\}}$.

As a consequence of these actions, we deduce the following.

\begin{propn}\label{infcores}
Suppose $(s\,|\,c)$ is an $l$-\mcd. Then $\mc sc$ is infinite \iff $s\neq1$.
\end{propn}

\begin{pf}
As noted above, when $s=1$ the only $(s\,|\,c)$-core is $\varnothing^l$. The case where $s\neq1$ follows from the discussion of actions above: the stabiliser of $\varnothing^l$ is easily seen to have infinite index in $\weyl s$ (in fact Hosaka \cite[Theorem 3.1]{hos} shows that a proper parabolic subgroup of any infinite irreducible Coxeter group has infinite index), so $\mc sc$ is in bijection with an infinite set.
\end{pf}

\section{Finiteness}\label{finitenesssec}

In this section we prove our main result: given a set $\calt$ of $l$-\mcds, we determine whether there are only finitely many multipartitions which are $(s\,|\,c)$-cores for all $(s\,|\,c)\in\calt$. We fix some notation.

\smallskip
\begin{mdframed}[innerleftmargin=3pt,innerrightmargin=3pt,innertopmargin=3pt,innerbottommargin=3pt,roundcorner=5pt,innermargin=-3pt,outermargin=-3pt]
\noindent\textbf{Notation in force for \cref{finitenesssec}:}
$l$ is a fixed positive integer, and $\calt$ is a set of $l$-\mcds. We write $\calt=\lset{(\sof st\,|\,\sof ct)}{t\in T}$ for an indexing set $T$.

We define $\smc\calt$ to be the intersection $\bigcap_{t\in T}\mc{\sof st}{\sof ct}$ (setting $\smc\calt=\calp^l$ when $\calt=\emptyset$), and we define $g(\calt)$ to be the \gcd of the integers in the set
\[
\lset{\sof st}{t\in T}\cup\lset{\sof ct_i-\sof ct_j-\sof cu_i+\sof cu_j}{t,u\in T,\ 1\ls i,j\ls l}.
\]
If the above set equals $\{0\}$ or is empty, then we set $g(\calt)=0$.
\end{mdframed}

\subsection{A simple criterion}

In this subsection we give a simple necessary condition for $\smc\calt$ to be finite. It will turn out that in almost all cases this condition is also sufficient. We begin with a useful lemma.

\begin{lemma}\label{divides}
Suppose $s,t\in\bbn\cup\{0\}$ and $c\in\bbz^l$, and that $s$ divides $t$. Then $\mc sc\subseteq\mc tc$.
\end{lemma}

Note that when we say $s$ divides $t$, we mean that $t=ns$ for some integer $n$, so we include the case $t=0$.

\begin{pf}
Since $s$ divides $t$, two nodes with the same $(t\,|\,c)$-residue must have the same $(s\,|\,c)$-residue. Hence two multipartitions with the same $(t\,|\,c)$-content have the same $(s\,|\,c)$-content. Now the result follows from the definition of $(s\,|\,c)$-cores.
\end{pf}

Now we can give our necessary condition for $\smc\calt$ to be finite.

\begin{cory}\label{gnot1}
Suppose $\smc\calt$ is finite. Then $g(\calt)=1$.
\end{cory}

\begin{pf}
Let $g=g(\calt)$, and observe that for any $t,u\in T$ there is $d\in\bbz$ such that we have $\sof ct_k\equiv\sof cu_k+d\ppmod g$ for all $k$. Hence by \cref{shiftmods}, $\mc g{\sof ct}=\mc g{\sof cu}$. In other words, the set $\mc g{\sof ct}$ is the same for every $t\in T$. By \cref{divides} $\mc g{\sof ct}\subseteq\mc{\sof st}{\sof ct}$, so $\smc\calt$ contains $\mc g{\sof ct}$. If $g\neq1$ then $\mc g{\sof ct}$ is infinite by \cref{infcores}, and hence so is $\smc\calt$.
\end{pf}

\subsection{The case where every $\sof st$ is zero}

In this subsection we assume that $\sof st=0$ for all $t\in T$. Perhaps surprisingly, this is the most complicated case.%, and we will use it in the next section where we address the case where $\sof st>0$ for some $t\in T$.

We begin with a simple construction of core multipartitions.

\begin{lemma}\label{xcore}
Suppose $c\in\bbz^l$, let $m=\max\{c_1,\dots,c_l\}$, and let $K=\lset{k\in\{1,\dots,l\}}{c_k=m}$. For any $n\in\bbn$ define a multipartition $\la$ by
\[
\la^{(k)}=
\begin{cases}
(n)&(k\in K)\\
\varnothing&(k\notin K).
\end{cases}
\]
Then $\la\in\mc0c$.
\end{lemma}

\begin{pf}
The $(0\,|\,c)$-content of $\la$ is $\{m^{|K|},(m+1)^{|K|},\dots,(m+n-1)^{|K|}\}$. Suppose $\mu$ is a multipartition with this $(0\,|\,c)$-content. Then $\mu^{(k)}=\varnothing$ for $k\notin K$, since $\mu$ has no nodes of residue less than $m$; for the same reason, $\mu^{(k)}_2=0$ for $k\in K$. Furthermore, $\mu^{(k)}_1\ls n$ for $k\in K$, because $\mu$ has no nodes of residue greater than $m+n-1$. The only possible $\mu$ satisfying these criteria is $\mu=\la$, so $\la$ is the unique multipartition with its $(0\,|\,c)$-content.
\end{pf}

Now we make a definition. Given $k\in\{1,\dots,l\}$, say that $k$ is
\begin{itemize}
\item
\emph{always maximal} if $\sof ct_k\gs\sof ct_j$ for all $t\in T$ and $j\in\{1,\dots,l\}$;
\item
\emph{sometimes maximal} if there is some $t\in T$ such that $\sof ct_k\gs\sof ct_j$ for all $j\in\{1,\dots,l\}$;
\item
\emph{never maximal} if for every $t\in T$ there is $j\in\{1,\dots,l\}$ with $\sof ct_k<\sof ct_j$.
\end{itemize}
We define \emph{always minimal}, \emph{never minimal} and \emph{sometimes minimal} similarly, with the inequalities reversed.

Say that $\calt$ satisfies \emph{condition X} if there is at least one $k\in\{1,\dots,l\}$ which is sometimes maximal but not always maximal, and at least one $k$ which is sometimes minimal but not always minimal.
\begin{comment}are $t,u,v,w\in T$ and $1\ls j,k\ls l$ such that
\begin{alignat*}4
\sof ct_j&\gs \sof ct_m&&\text{ for all }1\ls m\ls l,&\qquad \sof cu_j&<\sof cu_m&&\text{ for some }1\ls m\ls l,
\\
\sof cv_k&\ls \sof cv_m&&\text{ for all }1\ls m\ls l,&\qquad \sof cw_k&>\sof cw_m&&\text{ for some }1\ls m\ls l.
\end{alignat*}
Another way of saying this is: the set of elements of $\{1,\dots,l\}$ on which $\sof ct$ attains its maximum is not the same for all $t\in T$, and the set of elements of $\{1,\dots,l\}$ on which $\sof ct$ attains its minimum is not the same for all $t\in T$.
\end{comment}

Now we can state our main result for the case where every $\sof st$ equals $0$.

\smallskip
\begin{mdframed}[innerleftmargin=3pt,innerrightmargin=3pt,innertopmargin=-6pt,innerbottommargin=3pt,roundcorner=5pt,innermargin=-3pt,outermargin=-3pt]
\begin{thm}\label{mainall0}
Suppose $\calt=\lset{(0\,|\,\sof ct)}{t\in T}$ is a set of $l$-\mcds. Then $\smc\calt$ is finite \iff $g(\calt)=1$ and $\calt$ satisfies condition X.
\end{thm}
\end{mdframed}
\smallskip

\begin{eg}
Suppose $\calt=\{(0\,|\,1,3,0),(0\,|\,3,0,1)\}$. Then $g(\calt)=1$ and $\calt$ satisfies condition X. If $\la\in\smc\calt$, then by \cref{lev2} $(\la^{(1)},\la^{(2)})\in\mc0{1,3}\cap\mc0{3,0}$. \cref{00l2} below then tells us that $(\la^{(1)},\la^{(2)})\in\mc5{1,3}$, and in particular $\la^{(1)}$ and $\la^{(2)}$ are both $5$-cores. Similarly, $(\la^{(2)},\la^{(3)})\in\mc4{3,0}$, so $\la^{(2)}$ and $\la^{(3)}$ are both $4$-cores; since there are only finitely many $(4,5)$-cores, there are only finitely many possibilities for $\la^{(2)}$. It follows from \ref{sposbetacore} that for a given $5$-core $\la^{(2)}$ there are only finitely many bipartitions $(\la^{(1)},\la^{(2)})$ in $\mc5{1,3}$. So there are only finitely many possibilities for $\la^{(1)}$. Similarly, there are only finitely many possibilities for $\la^{(3)}$, and so $\smc\calt$ is finite.

In fact, we find that $\cardx{\smc\calt}=30$, with the largest tripartition in $\smc\calt$ being $((1^3),(3^2,1^3),(2^2))$.
\end{eg}

One direction of the proof is easy.

\begin{pf}[Proof of \cref{mainall0} (`only if' part)]
By \cref{gnot1} $\smc\calt$ is infinite if $g\neq1$. Now suppose $T$ does not satisfy condition X. This means either that every $k$ which is sometimes maximal is always maximal, or that every $k$ which is sometimes minimal is always minimal. We assume we are in the first case (the other case is similar). By \cref{xcore} the multipartition $\la$ given by
\[
\cpt\la k=
\begin{cases}
(n)&\text{if $k$ is always maximal}\\
\varnothing&\text{otherwise}
\end{cases}
\]
lies in $\smc\calt$ for every $n$, so $\smc\calt$ is infinite.
\end{pf}

Now we address the `if' part of \cref{mainall0}, which is considerably harder. Given $a>0$, let $\hk a\la$ denote the number of $a$-hooks of a partition $\la$. In particular, $\hk1\la$ is just the number of removable nodes of $\la$. The idea of the proof of \cref{mainall0} is to bound $\hk a{\cpt\la j}$ for $\la\in\smc\calt$, for each integer $a$ of the form $\left|\sof ct_j-\sof ct_k-\sof cu_j+\sof cu_k\right|$. The fact that these integers $a$ are coprime is then used to bound $\hk1{\cpt\la j}$. Condition X is then used to finish off the proof.

We start with a result on simultaneous core bipartitions which will also be useful in \cref{finitesec}.

\begin{lemma}\label{00l2}
Suppose $c_1,c_2,d_1,d_2\in\bbz$ with $c_1-c_2\gs0>d_1-d_2$, and let $a=c_1-c_2-d_1+d_2$. Then
\[
\mc0{c_1,c_2}\cap\mc0{d_1,d_2}=\mc a{c_1,c_2}.
\]
\end{lemma}

\begin{pf}
Since $c_1-c_2\equiv d_1-d_2\ppmod a$, we have $\mc a{c_1,c_2}=\mc a{d_1,d_2}$ by \cref{shiftmods}. Moreover, this set is contained in both $\mc0{c_1,c_2}$ and $\mc0{d_1,d_2}$ by \cref{divides}. So we just need to show that if $\la\in\mc0{c_1,c_2}\cap\mc0{d_1,d_2}$ then $\la\in\mc a{c_1,c_2}$. To see this, note that by \ref{s0betacore}
\[
\ber{\cpt\la1}{c_1}\supseteq\ber{\cpt\la2}{c_2},\qquad\ber{\cpt\la1}{d_1}\subseteq\ber{\cpt\la2}{d_2}
\]
so that
\[
\ber{\cpt\la1}{c_1-c_2}\supseteq\be{\cpt\la2}\supseteq\ber{\cpt\la1}{d_1-d_2}.
\]
The inequalities $c_1-c_2\gs0>d_1-d_2$ mean that the residue of $c_1-c_2$ modulo $a$ is $c_1-c_2$, so $\la\in\mc a{c_1,c_2}$ by \ref{sposbetacore}.
\end{pf}

We derive a simple consequence for simultaneous core multipartitions.

\begin{cory}\label{simbiopp2}
Suppose $\sof ct_j-\sof ct_k\gs0>\sof cu_j-\sof cu_k$ for some $t,u\in T$ and $1\ls j,k\ls l$. Let $a=\sof ct_j-\sof ct_k-\sof cu_j+\sof cu_k$. If $\la\in\smc\calt$, then $\cpt\la j$ and $\cpt\la k$ are $a$-cores.
\end{cory}

\begin{pf}
By \cref{lev2} the bipartition $(\cpt\la j,\cpt\la k)$ is both a $(0\,|\,\sof ct_j,\sof ct_k)$-core and an $(0\,|\,\sof cu_j,\sof cu_k)$-core. So by \cref{00l2} $(\cpt\la j,\cpt\la k)$ is an $(a\,|\,\sof cu_j,\sof cu_k)$-core, and in particular $\cpt\la j$ and $\cpt\la k$ are $a$-cores.
\end{pf}

Note that the difference in the signs of $c_1-c_2$ and $d_1-d_2$ is crucial in \cref{00l2}. In the absence of this hypothesis, the components of a bipartition in $\mc0{c_1,c_2}\cap\mc0{d_1,d_2}$ need not be $a$-cores. However, we can give a weaker result which shows that we can bound the number of $a$-hooks of each component. 

\begin{lemma}\label{simbibd1}
Suppose $c_1,c_2,d_1,d_2\in\bbz$ with $c_1-c_2>d_1-d_2\gs0$, and let $a=c_1-c_2-d_1+d_2$. If $\la\in\mc0{c_1,c_2}\cap\mc0{d_1,d_2}$, then $\hk a{\cpt\la k}\ls d_1-d_2$ for $k=1,2$.
\end{lemma}

\begin{pf}
We consider only $\cpt\la2$ (the proof for $\cpt\la1$ is similar)%, and as in the proof of \cref{simbiopp1} we assume $r=s=0$
. From \ref{s0betacore} and \cref{betaobs} we know that
\[
\ber{\cpt\la1}{c_1-c_2}=\be{\cpt\la2}\sqcup C,\qquad\ber{\cpt\la1}{d_1-d_2}=\be{\cpt\la2}\sqcup D
\]
for some sets $C,D$ of sizes $c_1-c_2,d_1-d_2$ respectively. Hence
\[
\rset{b-a}{b\in\be{\cpt\la2}\sqcup C}=\be{\cpt\la2}\sqcup D.
\]
This means that if $b\in\be{\cpt\la2}$ but $b-a\notin\be{\cpt\la2}$, then $b-a\in D$. Hence there are only $|D|=d_1-d_2$ possible values for $b$, so by \cref{corebeta} $\cpt\la2$ has at most $d_1-d_2$ $a$-hooks.
\end{pf}

Again, we note the consequences for multipartitions in $\calt$.

\begin{cory}\label{simbibd2}
Suppose $t,u\in T$ and $1\ls j,k\ls l$, and let $a=|\sof ct_j-\sof ct_k-\sof cu_j+\sof cu_k|$. If $a>0$ and $\la\in\smc\calt$, then
\[
\max\left\{\hk a{\cpt\la j},\hk a{\cpt\la k}\right\}\ls\min\left\{\left|\sof ct_j-\sof ct_k\right|,\left|\sof cu_j-\sof cu_k\right|\right\}.
\]
\end{cory}

\begin{pf}
By interchanging $j$ and $k$ or $t$ and $u$ if necessary, we can assume $\sof ct_j-\sof ct_k>\sof cu_j-\sof cu_k$ and $\sof ct_j-\sof ct_k\gs0$. If $\sof cu_j-\sof cu_k<0$, then the result follows from \cref{simbiopp2}, since then $\cpt\la j$ and $\cpt\la k$ are $a$-cores. So assume $\sof cu_j-\sof cu_k\gs0$. Since $(\cpt\la j,\cpt\la k)$ is both an $(\sof st\,|\,\sof ct_j,\sof ct_k)$-core and an $(\sof su\,|\,\sof cu_j,\sof cu_k)$-core, it is also a $(0\,|\,\sof ct_j,\sof ct_k)$-core and an $(0\,|\,\sof cu_j,\sof cu_k)$-core by \cref{divides}, so the result follows from \cref{simbibd1}.
\end{pf}

The preceding results show that for $\la\in\smc\calt$ the number of $a$-hooks of $\cpt\la j$ is bounded for each $a$ of the form $\left|\sof ct_j-\sof ct_k-\sof cu_j+\sof cu_k\right|$. We want to use this to show that $\hk1{\cpt\la j}$ is bounded. We do this via the following general result.
	
\begin{propn}\label{hookbd}
Suppose $P$ is a set of partitions, $A$ a set of coprime positive integers and $f:A\to\bbn$ a function such that $\hk a\la<f(a)$ for all $\la\in P$ and $a\in A$. Then there is $M\in\bbn$ such that $\hk1\la<M$ for all $\la\in P$.
\end{propn}

\begin{pf}
We assume that $A$ is finite; if it is not, we can certainly replace $A$ with a finite subset whose elements are still coprime. Since the elements of $A$ are coprime, we can find $G\in\bbn$ such that every integer greater than $G$ can be written as a sum of elements of $A$. Suppose for a contradiction that $\hk1\la$ is unbounded for $\la\in P$; then by \cref{corebeta} we can find, for any $M\in\bbn$, a partition $\la\in P$ and integers $b_1<\dots<b_M\in\be\la$ such that $b_1-1,\dots,b_M-1\notin\be\la$. Hence (letting $N=\inp{M/G}$) we can find $c_1<d_1<c_2<d_2<\dots<c_N<d_N$ such that for each $i$ we have $d_i-c_i>G$, $d_i\in\be\la$ and $c_i\notin\be\la$. But now by writing each $d_i-c_i$ as a sum of elements of $A$ and checking which integers between $c_i$ and $d_i$ lie in $\be\la$, we can find $c_i\ls e_i<f_i\ls d_i$ such that $f_i-e_i\in A$, $f_i\in\be\la$ and $e_i\notin\be\la$. Hence $\sum_{a\in A}\hk a\la\gs N$; taking $M$ such that $N>\sum_{a\in A}f(a)$ now gives a contradiction.
\end{pf}

%\begin{rmk}
%In fact, it is easy to prove a more general result: without the assumption that the elements of $A$ are coprime, let $g$ be their \gcd. Then one can show that the number of $g$-hooks of any partition in $P$ is bounded.
%\end{rmk}

As a consequence of this result, we see that when $g(\calt)=1$, the number of removable nodes of a multipartition in $\smc\calt$ is bounded, even without assuming Condition X. Now we use Condition X to complete the proof of the theorem. For this we need two more simple lemmas.

\begin{lemma}\label{fincorerem}
Suppose $s,b\in\bbn$. Then there are only finitely many $s$-core partitions having no more than $b$ removable nodes.
\end{lemma}

\begin{pf}
An $s$-core $\la$ satisfies $\la_i-\la_{i+1}<s$ for every $i$, since if $\la_i-\la_{i+1}\gs s$ then there is an $s$-hook contained in row $i$ of $[\la]$. So if $\la$ has no more than $b$ removable nodes, then $\la_1\ls(s-1)b$. Similarly, the length of the first column of $\la$ is at most $(s-1)b$, so $|\la|$ is bounded.
\end{pf}

\begin{lemma}\label{lengthmc}
If $\la\in\mc s{c_1,c_2}$ with $c_1\ls c_2$, then $\cpt\la1_1+c_1\ls\cpt\la2_1+c_2$.
\end{lemma}

\begin{pf}
By \cref{sbetacore,divides} we have $\ber{\cpt\la1}{c_1}\subseteq\ber{\cpt\la2}{c_2}$. In particular, $\cpt\la1_1+c_1\in\ber{\cpt\la2}{c_2}$, so there is $a\gs1$ such that $\cpt\la1_1-1+c_1=\cpt\la2_a-a+c_2$. But $\cpt\la2_a-a\ls\cpt\la2_1-1$, which gives the result. 
\end{pf}

Now we can proceed with the proof of the `if' part of \cref{mainall0}. Suppose $\sof st=0$ for all $t\in T$, and that $g(\calt)=1$ and $\calt$ satisfies condition X. Recall that $k\in\{1,\dots,l\}$ is \emph{\smx} if there is $t\in T$ such that $\sof ct_k\gs\sof ct_m$ for all $1\ls m\ls l$, and \emph{\smn} if there is $t\in T$ such that $\sof ct_k\ls\sof ct_m$ for all $1\ls m\ls l$.

\begin{lemma}\label{sometimesmx}
Suppose $k\in\{1,\dots,l\}$ is \smx or \smn. Then the set $\lset{\cpt\la k}{\la\in\smc\calt}$ is finite.
\end{lemma}

\begin{pf}
Given the assumption that $\sof st=0$ for every $t\in T$, $g(\calt)$ is the \gcd of the integers $\sof ct_k-\sof ct_j-\sof cu_k+\sof cu_j$ obtained as $j$ ranges over $\{1,\dots,l\}$ and $t,u$ range over $T$; so by assumption these integers are coprime. By \cref{simbibd2} if $\left|\sof ct_k-\sof ct_j-\sof cu_k+\sof cu_j\right|>0$ then the number of $\left|\sof ct_k-\sof ct_j-\sof cu_k+\sof cu_j\right|$-hooks of $\la$ is bounded as $\la$ ranges over $\smc\calt$. So if we let
\[
A=\lset{\left|\sof ct_k-\sof ct_j-\sof cu_k+\sof cu_j\right|}{t,u\in T,\ j\in\{1,\dots,l\}}\setminus\{0\}
\]
and $P=\lset{\cpt\la k}{\la\in\smc\calt}$, then $A$ and $P$ satisfy the hypotheses of \cref{hookbd}. So the number of $1$-hooks (i.e.\ the number of removable nodes) of a partition in $P$ is bounded, by $b$ say.

Now Condition X together with the fact that $k$ is \smx or \smn implies that there are $t,u\in T$ and $j\in\{1,\dots,l\}$ such that either $\sof ct_k-\sof ct_j\gs0>\sof cu_k-\sof cu_j$ or $\sof ct_j-\sof ct_k\gs0>\sof cu_j-\sof cu_k$. If we let $a=\left|\sof ct_k-\sof ct_j-\sof cu_k+\sof cu_j\right|$, then by \cref{simbiopp2} $\cpt\la k$ is an $a$-core for every $\la\in\smc\calt$. Since $a>0$, \cref{fincorerem} gives the result.
\end{pf}

Now we can complete the proof.

\begin{pf}[Proof of \cref{mainall0} (`if' part)]
Suppose $g(\calt)=1$ and $\calt$ satisfies condition X. To show that $\smc\calt$ is finite, it suffices to show that for every $k\in\{1,\dots,l\}$ the set $\lset{\cpt\la k}{\la\in\smc\calt}$ is finite. We have proved this when $k$ is \smx or \smn, so assume $k$ is never maximal and never minimal. The fact that $k$ is never maximal means that there is $j$ which is \smx and $t\in T$ such that $\sof ct_k\ls \sof ct_j$. If $\la\in\smc\calt$, then $(\cpt\la k,\cpt\la j)\in\mc0{\sof ct_k,\sof ct_j}$, so by \cref{lengthmc} $\cpt\la k_1\ls\cpt\la j_1+\sof ct_j-\sof ct_k$. Since by \cref{sometimesmx} there are only finitely many possible $\cpt\la j$, this means that $\cpt\la k_1$ is bounded as $\la$ ranges over $\smc\calt$. Similarly (using the fact that $k$ is never minimal) the first column of $\cpt\la k$ is bounded, so there are only finitely many possible~$\cpt\la k$.
\end{pf}

\subsection{The case where $\sof st>0$ for some $t$}\label{finitesec}

In this subsection we complete the analysis of when $\smc\calt$ is finite by considering the case where $\sof st>0$ for some $t\in T$. The statement here is simpler.

\smallskip
\begin{mdframed}[innerleftmargin=3pt,innerrightmargin=3pt,innertopmargin=-6pt,innerbottommargin=3pt,roundcorner=5pt,innermargin=-3pt,outermargin=-3pt]
\begin{thm}\label{mainnotall0}
Suppose $\calt=\lset{(\sof st\,|\,\sof ct)}{t\in T}$ is a set of $l$-\mcds with $\sof st>0$ for at least one $t\in T$. Then $\smc\calt$ is finite \iff $g(\calt)=1$.
\end{thm}
\end{mdframed}
\smallskip

We can deduce \cref{mainnotall0} fairly easily from \cref{mainall0}. To begin with, we use \cref{00l2} to express $\mc sc$ for any $s,c$ as an intersection of sets $\mc0d$.

\begin{propn}\label{intersect}
Suppose $(s\,|\,c)$ is an $l$-\mcd with $s>0$. Then
\[
\mc sc=\bigcap_{d\in\bbz^l}\mc0{c+sd}.
\]
\end{propn}

\begin{pf}
For each $d\in\bbz^l$ we have $\mc sc=\mc s{c+sd}\subseteq\mc0{c+sd}$ by \cref{shiftmods,divides}, so the left-hand side is contained in the right-hand side. For the opposite inclusion, suppose $\la\in\mc0{c+sd}$ for every $d\in\bbz^l$. Given $1\ls j<k\ls l$, we can find $d,e\in\bbz^l$ such that $0\ls (c_j+sd_j)-(c_k+sd_k)<s$ and $e_j-e_k=d_j-d_k-1$. Then $(\cpt\la j,\cpt\la k)\in\mc0{c_j+sd_j,c_k+sd_k}\cap\mc0{c_j+se_j,c_k+se_k}=\mc s{c_j,c_k}$, by \cref{lev2,00l2}. Since this is true for every $j,k$, we have $\la\in\mc sc$ by \cref{lev2}.
\end{pf}

\begin{rmk}
In fact, one can write $\mc sc=\bigcap_{d\in M}\mc0{c+sd}$ for a much smaller subset $M$ of $\bbz^l$: it is possible to take $|M|=l$. But it is easier for us to take $M$ to be the whole of $\bbz^l$ as in \cref{intersect}.
\end{rmk}

This yields the following.% \lcnamecref{change0}, which is enough to deduce \cref{mainnotall0} from \cref{mainall0}.

\begin{propn}\label{change0}
Suppose $\sof st>0$ for at least one $t\in T$. Then there is a set $\calu=\lset{(0\,|\,\sof cu)}{u\in U}$ of $l$-\mcds such that:
\begin{enumerate}
\item\label{hatx}
$\calu$ satisfies condition X;
\item\label{ghatg}
$g(\calu)=g(\calt)$;
\item\label{smchatsmc}
$\smc\calu=\smc\calt$.
\end{enumerate}
\end{propn}

\begin{pf}
%Let $T_0=\rset{t\in T}{\sof st=0}$ and $T_1=T\setminus T_0$. %For each $t\in T_1$, let $\sof Mt=\sof ct+\sof st\bbz^l$. 
%Now define
Define
\[
\calu=\lset{(0\,|\,\sof ct+\sof std)}{t\in T,\ d\in\bbz^l}.
\]
Now we check the conditions in the \lcnamecref{change0}.
\begin{enumerate}
\item
By assumption there is $t\in T$ such that $\sof st>0$. For any $1\ls j<k\ls l$ we can easily find $d,e\in\bbz^l$ such that $\sof ct_j+\sof std_j>\sof ct_k+\sof std_k$ and $\sof ct_j+\sof ste_j<\sof ct_k+\sof ste_k$. This shows that no $k\in\{1,\dots,l\}$ is always maximal or always minimal for $\calu$, which \textit{a fortiori} gives condition X for $\calu$.
\item
By definition $g(\calt)$ is the \gcd of the integers in the set
\[
\lset{\sof st}{t\in T}\cap\lset{\sof ct_i-\sof ct_j-\sof cu_i+\sof cu_j}{t,u\in T,\ 1\ls i,j\ls l},
\]
while $g(\calu)$ is the \gcd of the integers in the set
\[
\lset{\sof ct_i-\sof ct_j-\sof cu_i+\sof cu_j+a\sof st+b\sof su}{t,u\in T,\ 1\ls i,j\ls l,\ a,b\in\bbz}.
\]
It is easy to see that these \gcds are the same.
\item
This follows from \cref{intersect}.\qedhere
\end{enumerate}
\end{pf}

\begin{pf}[Proof of \cref{mainnotall0}]
The `only if' part is \cref{gnot1}. For the `if' part, suppose $g(\calt)=1$, and let $\calu$ be as in \cref{change0}. Then by \cref{mainall0} $\smc\calt=\smc\calu$ is finite.
\end{pf}

\section{Enumeration of simultaneous core multipartitions}

An early success in the study of simultaneous core partitions was Anderson's Theorem \cite[Theorems 1 \& 3]{and} that when $s$ and $t$ are coprime, the number of $(s,t)$-cores is the \emph{rational Catalan number} $\mfrac1{s+t}\mbinom{s+t}s$. Extending this to the enumeration of partitions in $\cores{s_1}\cap\dots\cap\cores{s_r}$ for coprime integers $s_1,\dots,s_r$ with $r\gs3$ seems to be much more difficult, although various special cases have been addressed in the recent literature \cite{hana,amd,amle,wang,xiong}.

Naturally, one can extend these enumerative questions to simultaneous core multipartitions: in particular, given a set $\calt$ of $l$-\mcds such that $\smc\calt$ is finite (as determined by \cref{mainall0,mainnotall0}), what is $\cardx{\smc\calt}$? This question seems to be very hard to answer in general; the proofs of \cref{mainall0,mainnotall0} do not give anything like an efficient algorithm for calculating $\smc\calt$, so it is difficult even to gather data. In this section we address the very simplest case, where $l=\card\calt=2$. Even here the enumeration question is difficult to answer, and we restrict to two particular subcases.

If $l=\cardx\calt=2$, we can assume (in view of \cref{shiftmods}) that
\[
\calt=\{(s\,|\,0,a),(t\,|\,0,b)\}
\]
with $s,t\in\bbn\cup\{0\}$ and $a,b\in\bbz$. Moreover, if $s>0$ then we can take $0\ls a<s$, and similarly for $t$ and $b$.

\subsection{The case where $s$ divides $t$}\label{divsec}

In this subsection we take $\calt$ as above with $s$ dividing $t$. We start with the case $s=t=0$. In this case $g(\calt)=|a-b|$, so we need $|a-b|=1$ in order to have $\smc\calt$ finite. But we also need $\calt$ to satisfy condition X, which means that $a$ or $b$ equals $0$. Now we have the following result.

\begin{propn}\label{s0t0}
Suppose $\calt=\{(0\,|\,0,a),(0\,|\,0,b)\}$, with $\{|a|,|b|\}=\{0,1\}$. Then $\cardx{\smc\calt}=1$.
\end{propn}

\begin{pf}
We assume $a=0$ and $b=1$ (the other cases follow symmetrically). Suppose $(\la,\mu)$ is a bipartition lying in $\smc\calt$; we will show that $\la=\mu=\varnothing$. By \ref{s0betacore} the fact that $(\la,\mu)\in\mc0{0,0}$ says that $\be\la=\be\mu$; since a partition can be recovered from its beta-set, we obtain $\la=\mu$. Now the fact that $(\la,\la)\in\mc0{0,1}$ gives $\be\la\subseteq\ber\la1$; by \cref{betaobs} this means that $\ber\la1=\be\la\cup\{b\}$ for some integer $b$. In fact it is easy to see that $b$ must equal $\la_1$ (since this lies in $\ber\la1$ and is larger than the largest element $\la_1-1$ of $\be\la$). Hence we have $\ber\la1\setminus\{\la_1\}=\be\la$; writing the elements of these sets in decreasing order, we obtain
\[
\la_2-1=\la_1-1,\qquad\la_3-2=\la_2-2,\qquad\la_4-3=\la_3-3,\dots
\]
so that $\la_1=\la_2=\la_3=\dots$, and therefore $\la=\varnothing$.
\end{pf}

Now we consider the case where $s,t>0$. We will deduce our main result here as a special case of a more general result. So to begin with we do not assume that $s$ divides $t$, and we let $g$ be the \gcd of $s$ and $t$ throughout this section. We will restrict attention to bipartitions $(\la,\mu)$ for which both $\la,\mu$ are $g$-cores. Let $\cores g^2$ denote the set of such bipartitions.

Let $U_g^{s,a}$ denote the set of all tuples $u=\ltup{u_i}{i\in\zgz}$ of integers with $\sum_i u_i=a$ and $0\ls u_i\ls s/g$ for each $i\in\zgz$. By a simple application of the Inclusion--Exclusion Principle,
\[
\card{U_g^{s,a}}=\sum_{d\gs0}(-1)^d\binom gd\binom{a+g-1-d(1+s/g)}{g-1}.
\]
Now we can state our main theorem in this section.

\begin{thm}\label{countss}
Suppose $0\ls a<s$ and $0\ls b<t$. Let $g=\ngcd(s,t)$, and assume $g$ and $a-b$ are coprime. Then
\[
\cardx{\bc sa\cap\bc tb\cap\smash{\cores g^2}}=\frac1g\card{U_g^{s,a}}\card{U_g^{t,b}}.
\]
In particular, if $s$ divides $t$, then
\[
\cardx{\bc sa\cap\bc tb}=\frac1s\binom sa\card{U_s^{t,b}}.
\]
\end{thm}

We remark that in the very special case where $s=t$, we get the even simpler formula
\[
\cardx{\bc sa\cap\bc sb}=\frac1s\binom sa\binom sb.
\]

To prove \cref{countss}, we use a slightly different version of \ref{sposbetacore} to characterise core bipartitions. Suppose $\la$ is an $s$-core. For each $i\in\zsz$, let $\se i\la$ be the smallest element of $i$ not contained in $\be\la$. The set $\se s\la=\lset{\se i\la}{i\in\zsz}$ is referred to as the \emph{$s$-set} of $\la$; these sets were studied extensively in \cite{mfcores,mfgencores,mfwtarmstrong}. Observe that $\se s\la$ is a set of $s$ integers which are pairwise incongruent modulo $s$ and sum to $\binom s2$. Conversely, any such set of integers is the $s$-set of a unique $s$-core.

The following lemma, which follows easily from the definition, shows how to obtain the $s$-set of a $g$-core from its $g$-set.

\begin{lemma}\label{gsetsset}
Suppose $s,g$ are integers with $g\mid s$, and $\la\in\cores g$. Then
\[
\se s\la=\lset{\se i\la+kg}{i\in\zgz,\ 0\ls k<s/g}.
\]
\end{lemma}

Using $s$-sets, we can give a different version of \ref{sposbetacore} (in fact, this is much closer to the original version of this result in \cite{mfwt}).

\begin{propn}\label{sposbetacore2}
Suppose $\la,\mu\in\calp$ and $0\ls a<s$. Then  $(\la,\mu)\in\mc s{0,a}$ \iff $\la,\mu\in\cores s$ and
\[
\se i\mu+a\in\{\se{i+a}\la,\se{i+a}\la+s\}
\]
for each~$i\in\zsz$.
\end{propn}

\begin{pf}
This follows easily from \ref{sposbetacore}.
\end{pf}

In order to use \cref{sposbetacore2} to prove \cref{countss}, we want to consider bipartitions $(\la,\mu)\in\bc sa\cap \cores g^2$. So suppose $\la,\mu\in\cores g$. \cref{sposbetacore2} says that $(\la,\mu)\in\bc sa$ \iff $\se i\mu+a-\se{i+a}\la$ equals either $s$ or $0$ for each $i$. Since $g\mid s$ and $\la,\mu\in\cores s$, we have
\begin{align*}
\se s\la&=\lset{\se{i+a}\la+kg}{i\in\zgz,\ 0\ls k<s/g},\\\se s\mu&=\lset{\se i\mu+kg}{i\in\zgz,\ 0\ls k<s/g}
\end{align*}
by \cref{gsetsset}. So if $(\la,\mu)\in\bc sa$, then for each $i\in\zgz$ there is an integer $u_i\in\{0,\dots,s/g\}$ such that $\se i\mu+a=\se{i+a}\la+gu_i$. Summing over $i$ and using the fact that $\sum_{i\in\zgz}\se i\la=\sum_{i\in\zgz}\se i\mu$, we find that $\sum_{i\in\zgz}u_i=a$, so that the tuple $u=\lset{u_i}{i\in\zgz}$ lies in $u\in U_g^{s,a}$.

So if we define $\sigma(\la,\mu)=u$, we obtain a function
\begin{align*}
\sigma:\bc sa\cap\cores g^2\longrightarrow U_g^{s,a}.
\\
\intertext{Doing the same with $s,a$ replaced by $t,b$, we get another function}
\tau:\bc tb\cap\cores g^2\longrightarrow U_g^{t,b}.
\end{align*}

\begin{pf}[Proof of \cref{countss}]
We want to consider the images of the maps $\sigma,\tau$ defined above, and for this we need some more notation. For any tuple $u=\ltup{u_i}{i\in\zgz}$ and any $c\in\zgz$ define the tuple $u(+c)$ by $u(+c)_i=u_{i+c}$ for each $i\in\zgz$.

To prove the theorem we will prove the following claim: given $u\in U_g^{s,a}$ and $v\in U_g^{t,b}$, there is a unique bipartition $(\la,\mu)\in\mc s{0,a}\cap\mc t{0,b}\cap\cores g^2$ and a unique $c\in\zsz$ such that $\sigma(\la,\mu)=u(+c)$ and $\tau(\la,\mu)=v(+c)$.

Our first aim is to find integers $x_i,y_i$ for $i\in\zgz$ such that
\[
y_i=x_{i+a}-a+gu_i=x_{i+b}-b+gv_i\tag*{(\textasteriskcentered)}
\]
for all $i$. In fact, this is straightforward: we just fix $k\in\bbz$ and set
\begin{align*}
x_{d(b-a)+g\bbz}&=k+d(b-a)+g\sum_{j=0}^{d-1}u_{j(b-a)+g\bbz}-g\sum_{j=1}^dv_{j(b-a)+g\bbz}\\
y_{d(b-a)-a+g\bbz}&=k+d(b-a)-a+g\sum_{j=0}^du_{j(b-a)+g\bbz}-g\sum_{j=1}^dv_{j(b-a)+g\bbz}
\end{align*}
for all $0\ls d<g$. Since $b-a$ and $g$ are coprime, this uniquely defines $x_i$ and $y_i$ for every $i\in\zgz$, and it is easy to see that (\textasteriskcentered) is satisfied. Moreover, apart from the choice of $k$, these are the unique integers $x_i,y_i$ satisfying (\textasteriskcentered): once $x_{0+g\bbz}=k$ is chosen, (\textasteriskcentered) forces the choice of $y_{-a+g\bbz},x_{b-a+g\bbz},y_{b-2a+g\bbz},x_{2b-2a+g\bbz},\dots$, so that $x_i$ and $y_i$ are forced for every $i$.

Now observe that the integers $x_i$ are pairwise incongruent modulo $g$, so in particular sum to $\binom g2$ modulo $g$. Changing $k$ by $1$ changes this sum by $g$, and therefore there is a unique choice of $k$ (which we fix henceforth) such that $\sum_ix_i=\binom g2$. This also gives $\sum_iy_i=\binom g2$, so $\lset{x_i}{i\in\zgz}$ and $\lset{y_i}{i\in\zgz}$ are the $g$-sets of $g$-cores $\la$ and $\mu$ respectively. Since $x_i\in i+k$, we have $\se i\la=x_{i-k}$, and similarly $\se i\mu=y_{i-k}$, for each $i\in\zgz$, and hence
\[
\se i\mu+a-\se{i+a}\la=y_{i-k}-x_{i+a-k}+a=gu_{i+k},
\]
so that $(\la,\mu)\in\mc s{0,a}$ with $\sigma(\la,\mu)=u(+c)$, where $c=k+g\bbz$. Similarly $(\la,\mu)\in\bc tb$ with $\tau(\la,\mu)=v(+c)$, so we have the required $\la,\mu,c$. Moreover, the integers $x_i,y_i$ can be recovered from $\la,\mu,c$, so (by the statement above about the uniqueness of $x_i,y_i$) we have uniqueness for $\la,\mu,c$.

As a consequence of this claim, we find that $\cardx{\bc sa\cap\bc tb\cap\smash{\cores g^2}}$ equals $\frac1g$ times the number of choices of $u,v$. $u$ can be chosen in $\card{U_g^{s,a}}$ ways, and $v$ in $\card{U_g^{t,b}}$ ways, giving the result.

For the special case where $s$ divides $t$, we have $g=s$, so that $\card{U_g^{s,a}}=\binom sa$. Furthermore, $\bc sa\subseteq\cores g^2$ by \cref{multicore}, and the result follows.
\end{pf}

\begin{eg}
Take $s=3$, $t=9$, $a=1$ and $b=5$. The twelve bipartitions $(\la,\mu)\in\bc31\cap\bc95$ are given by the following table, where we give $\se3\la$, $\se3\mu$, $\sigma(\la,\mu)$, $\tau(\la,\mu)$, writing each $\bbz/3\bbz$-tuple $u$ in the form $(u_{0+3\bbz},u_{1+3\bbz},u_{2+3\bbz})$. We see that up to simultaneous cyclic permutation, each pair in $U_3^{3,1}\times U_3^{9,5}$ occurs once as $(\sigma(\la,\mu),\tau(\la,\mu))$.
\[
\begin{array}{cccccc}\hline
\la&\mu&\se3\la&\se3\mu&\sigma(\la,\mu)&\tau(\la,\mu)\\\hline
\varnothing&\varnothing&\{0,1,2\}&\{0,1,2\}&(0,0,1)&(1,2,2)\\
\varnothing&(1)&\{0,1,2\}&\{3,1,-1\}&(1,0,0)&(2,2,1)\\
(1)&\varnothing&\{3,1,-1\}&\{0,1,2\}&(0,1,0)&(2,1,2)\\
\varnothing&(2)&\{0,1,2\}&\{0,4,-1\}&(0,1,0)&(1,3,1)\\
(1^2)&\varnothing&\{3,-2,2\}&\{0,1,2\}&(1,0,0)&(1,1,3)\\
(1)&(1^2)&\{3,1,-1\}&\{3,-2,2\}&(1,0,0)&(3,0,2)\\
(2)&(1)&\{0,4,-1\}&\{3,1,-1\}&(0,1,0)&(3,2,0)\\
(2)&(1^2)&\{0,4,-1\}&\{3,-2,2\}&(0,0,1)&(3,1,1)\\
(1)&(3,1)&\{3,1,-1\}&\{0,-2,5\}&(0,0,1)&(2,0,3)\\
(2,1^2)&(1)&\{-3,4,2\}&\{3,1,-1\}&(0,0,1)&(2,3,0)\\
(1^2)&(2,1^2)&\{3,-2,2\}&\{-3,4,2\}&(0,1,0)&(0,2,3)\\
(3,1)&(2)&\{0,-2,5\}&\{0,4,-1\}&(1,0,0)&(0,3,2)\\\hline
\end{array}
\]
\end{eg}

\begin{rmk}
To complete the study of the situation where $s$ divides $t$, it remains to consider the case where $s>0$ and $t=0$. We deal with this case as a limiting case of \cref{countss}. So take $s,t,a,b$, with $t=ns$ for $n\in\bbn$. For large $n$ (in fact, for $n\gs b$), the value of $\card{U_s^{ns,b}}$ stabilises at $\binom{b+s-1}{s-1}$. In addition, one can see from the proof of \cref{countss} that the set $\bc sa\cap\bc{ns}b$ stabilises; call this limiting set $\calc$. We claim that $\bc sa\cap\bc0b=\calc$. By \cref{divides}, $\bc{ns}b\subseteq\bc0b$ for every $n$, so we have $\calc\subseteq\bc sa\cap\bc0b$. On the other hand, given a bipartition $(\la,\mu)$ and given $N$ sufficiently large relative to $(\la,\mu)$, we have $(\la,\mu)\in\bc Nb$ \iff $(\la,\mu)\in\bc0b$: we just take $N$ large enough that any two nodes which can occur as nodes of bipartitions of size $|\la|+|\mu|$ and which have the same $(N\,|\,0,b)$-residue must also have the same $(0\,|\,0,b)$-residue. So if $(\la,\mu)\notin\calc$, then $(\la,\mu)\notin\bc sa\cap\bc{ns}b$ for sufficiently large $n$, so that $(\la,\mu)\notin\bc sa\cap\bc0b$. Hence $\bc sa\cap\bc0b\subseteq\calc$, so $\bc sa\cap\bc0b=\calc$ as required.

So we deduce that
\[
\cardx{\bc sa\cap\bc 0b}=\frac1s\binom sa\binom{b+s-1}{s-1}.
\]
\end{rmk}

\subsection{The case $0\ls a=b<s,t$}\label{aasec}

Now we consider the case where the residue of $a$ modulo $s$ is the same as the residue of $b$ modulo $t$. In this case, we may assume that $0\ls a=b<s,t$.

\begin{thm}\label{countaa}
Suppose $0\ls a<s\ls t$, and that $s$ and $t$ are coprime. Then
\[
\cardx{\bc sa\cap\bc ta}=\frac{(s+t-a-1)!}{a!(s-a)!(t-a)!}.
\]
\end{thm}

In order to prove \cref{countaa}, we recall the \emph{$(s,t)$-lattice} used in Anderson's proof of her theorem. This is a diagram of $\bbz^2$, with the point $(x,y)$ replaced by the integer $sx+ty$. For example, part of the $(3,5)$-lattice is drawn as follows.
\[
\begin{tikzpicture}[scale=.9]
\draw(-3,-1)node{$-30$};
\draw(-3,0)node{$-25$};
\draw(-3,1)node{$-20$};
\draw(-3,2)node{$-15$};
\draw(-3,3)node{$-10$};
\draw(-3,4)node{$-5$};
\draw(-3,5)node{$0$};
\draw(-2,-1)node{$-27$};
\draw(-2,0)node{$-22$};
\draw(-2,1)node{$-17$};
\draw(-2,2)node{$-12$};
\draw(-2,3)node{$-7$};
\draw(-2,4)node{$-2$};
\draw(-2,5)node{$3$};
\draw(-1,0)node{$-19$};
\draw(-1,1)node{$-14$};
\draw(-1,2)node{$-9$};
\draw(-1,3)node{$-4$};
\draw(-1,4)node{$1$};
\draw(-1,5)node{$6$};
\draw(0,0)node{$-16$};
\draw(0,1)node{$-11$};
\draw(0,2)node{$-6$};
\draw(0,3)node{$-1$};
\draw(0,4)node{$4$};
\draw(0,5)node{$9$};
\draw(1,0)node{$-13$};
\draw(1,1)node{$-8$};
\draw(1,2)node{$-3$};
\draw(1,3)node{$2$};
\draw(1,4)node{$7$};
\draw(1,5)node{$12$};
\draw(2,0)node{$-10$};
\draw(2,1)node{$-5$};
\draw(2,2)node{$0$};
\draw(2,3)node{$5$};
\draw(2,4)node{$10$};
\draw(2,5)node{$15$};
\draw(3,0)node{$-7$};
\draw(3,1)node{$-2$};
\draw(3,2)node{$3$};
\draw(3,3)node{$8$};
\draw(3,4)node{$13$};
\draw(3,5)node{$18$};
\draw(4,0)node{$-4$};
\draw(4,1)node{$1$};
\draw(4,2)node{$6$};
\draw(4,3)node{$11$};
\draw(4,4)node{$16$};
\draw(4,5)node{$21$};
\draw(5,0)node{$-1$};
\draw(5,1)node{$4$};
\draw(5,2)node{$9$};
\draw(5,3)node{$14$};
\draw(5,4)node{$19$};
\draw(5,5)node{$24$};
%\draw(5,6)node{$29$};
\draw(6,0)node{$2$};
\draw(6,1)node{$7$};
\draw(6,2)node{$12$};
\draw(6,3)node{$17$};
\draw(6,4)node{$22$};
\draw(6,5)node{$27$};
%\draw(6,6)node{$32$};
\draw(7,0)node{$5$};
\draw(7,1)node{$10$};
\draw(7,2)node{$15$};
\draw(7,3)node{$20$};
\draw(7,4)node{$25$};
\draw(7,5)node{$30$};
%\draw(7,6)node{$35$};
\draw(-1,-1)node{$-24$};
\draw(0,-1)node{$-21$};
\draw(1,-1)node{$-18$};
\draw(2,-1)node{$-15$};
\draw(3,-1)node{$-12$};
\draw(4,-1)node{$-9$};
\draw(5,-1)node{$-6$};
\draw(6,-1)node{$-3$};
\draw(7,-1)node{$0$};
\end{tikzpicture}
\]
Note that the $(s,t)$-lattice is periodic: it is unchanged under translations by multiples of the vector $(t,-s)$. To construct the \emph{$(s,t)$-diagram} (sometimes called the $(s,t)$-abacus diagram) of a partition $\la$, one simply colours or circles the integers lying in $\be\la$. By \cref{corebeta}, the condition that $\la$ is an $(s,t)$-core is then simply that each coloured position has coloured positions both below and to the left. Part of the $(3,5)$-diagram of the $(3,5)$-core $(1)$ is as follows.
\[
\begin{tikzpicture}[scale=.9]
\path[fill=gray!50!white](-3.5,5.5)--++(1,0)--++(0,-1)--++(1,0)--++(0,-1)--++(1,0)--++(0,-1)--++(3,0)--++(0,-1)--++(1,0)--++(0,-1)--++(1,0)--++(0,-1)--++(3,0)--++(0,-1)--++(-11,0);
\draw(-3,-1)node{$-30$};
\draw(-3,0)node{$-25$};
\draw(-3,1)node{$-20$};
\draw(-3,2)node{$-15$};
\draw(-3,3)node{$-10$};
\draw(-3,4)node{$-5$};
\draw(-3,5)node{$0$};
\draw(-2,-1)node{$-27$};
\draw(-2,0)node{$-22$};
\draw(-2,1)node{$-17$};
\draw(-2,2)node{$-12$};
\draw(-2,3)node{$-7$};
\draw(-2,4)node{$-2$};
\draw(-2,5)node{$3$};
\draw(-1,0)node{$-19$};
\draw(-1,1)node{$-14$};
\draw(-1,2)node{$-9$};
\draw(-1,3)node{$-4$};
\draw(-1,4)node{$1$};
\draw(-1,5)node{$6$};
%\draw(-1,6)node{$11$};
\draw(0,0)node{$-16$};
\draw(0,1)node{$-11$};
\draw(0,2)node{$-6$};
\draw(0,3)node{$-1$};
\draw(0,4)node{$4$};
\draw(0,5)node{$9$};
%\draw(0,6)node{$14$};
\draw(1,0)node{$-13$};
\draw(1,1)node{$-8$};
\draw(1,2)node{$-3$};
\draw(1,3)node{$2$};
\draw(1,4)node{$7$};
\draw(1,5)node{$12$};
%\draw(1,6)node{$17$};
\draw(2,0)node{$-10$};
\draw(2,1)node{$-5$};
\draw(2,2)node{$0$};
\draw(2,3)node{$5$};
\draw(2,4)node{$10$};
\draw(2,5)node{$15$};
%\draw(2,6)node{$20$};
\draw(3,0)node{$-7$};
\draw(3,1)node{$-2$};
\draw(3,2)node{$3$};
\draw(3,3)node{$8$};
\draw(3,4)node{$13$};
\draw(3,5)node{$18$};
%\draw(3,6)node{$23$};
\draw(4,0)node{$-4$};
\draw(4,1)node{$1$};
\draw(4,2)node{$6$};
\draw(4,3)node{$11$};
\draw(4,4)node{$16$};
\draw(4,5)node{$21$};
%\draw(4,6)node{$26$};
\draw(5,0)node{$-1$};
\draw(5,1)node{$4$};
\draw(5,2)node{$9$};
\draw(5,3)node{$14$};
\draw(5,4)node{$19$};
\draw(5,5)node{$24$};
%\draw(5,6)node{$29$};
\draw(6,0)node{$2$};
\draw(6,1)node{$7$};
\draw(6,2)node{$12$};
\draw(6,3)node{$17$};
\draw(6,4)node{$22$};
\draw(6,5)node{$27$};
\draw(7,0)node{$5$};
\draw(7,1)node{$10$};
\draw(7,2)node{$15$};
\draw(7,3)node{$20$};
\draw(7,4)node{$25$};
\draw(7,5)node{$30$};
\draw(-1,-1)node{$-24$};
\draw(0,-1)node{$-21$};
\draw(1,-1)node{$-18$};
\draw(2,-1)node{$-15$};
\draw(3,-1)node{$-12$};
\draw(4,-1)node{$-9$};
\draw(5,-1)node{$-6$};
\draw(6,-1)node{$-3$};
\draw(7,-1)node{$0$};
\draw[very thick](-3.5,5.5)--++(1,0)--++(0,-1)--++(1,0)--++(0,-1)--++(1,0)--++(0,-1)--++(3,0)--++(0,-1)--++(1,0)--++(0,-1)--++(1,0)--++(0,-1)--++(3,0)--++(0,-1);
\end{tikzpicture}
\]
Now consider the boundary between the coloured and uncoloured parts of the diagram. The condition that $\la$ is an $(s,t)$-core means that this path consists only of steps to the right and steps down. Moreover, it is periodic, with each period consisting of $t$ steps to the right and $s$ steps down. We can encode this boundary path by writing down one period; of course, any cyclic permutation of this period will encode the same periodic boundary path. For example, we can encode the boundary path in the diagram above by (any cyclic permutation of) the sequence DRDRDRRR.

Conversely, any cyclic sequence comprising $s$ Rs and $t$ Ds yields the $(s,t)$-diagram of an $(s,t)$-core: if we draw the corresponding periodic path in the $(s,t)$-lattice, then the set of integers below and to the left of the path is the shifted beta-set of an $(s,t)$-core. Translating the path to a different position just changes the shift of the beta-set, without changing the partition.

As a consequence, we find that the number of $(s,t)$-cores equals the number of arrangements of $s$ Rs and $t$ Ds modulo cyclic shifts, which yields Anderson's Theorem.

Now we extend these ideas to the setting of \cref{countaa}. Suppose we have $0\ls a<s,t$, and that $(\la,\mu)\in\bc sa\cap\bc ta$. Consider the shifted beta-set $\ber\mu a$. By \ref{sposbetacore}, this is obtained from $\be\la$ by adding $a$ integers $x_1,\dots,x_a$, with $x_j-s,x_j-t\in\be\la$ for each $j$. So drawing the $a$-shifted $(s,t)$-diagram of $\mu$ (i.e.\ colouring the elements of $\ber\mu a$) amounts to taking the $(s,t)$-diagram of $\la$ and additionally colouring $a$ integers each of which has coloured integers both immediately below and immediately to the left. For example, take $(s,t,a)=(3,5,2)$, and $(\la,\mu)=((1),(2))$. Combining the $(3,5)$-diagram of $\la$ and the $2$-shifted $(3,5)$-diagram of $\mu$, we get the following picture (in which we use a lighter colour for the additional positions coloured in $\mu$).
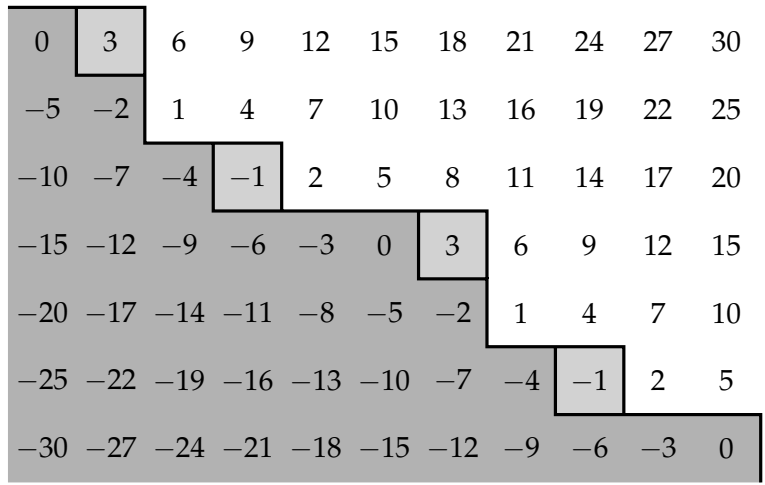
\[
\begin{tikzpicture}[scale=.9]
\path[fill=gray!35!white](-.5,2.5)--++(1,0)--++(0,1)--++(-1,0);
\path[fill=gray!35!white](4.5,-.5)--++(1,0)--++(0,1)--++(-1,0);
\path[fill=gray!35!white](2.5,1.5)--++(1,0)--++(0,1)--++(-1,0);
\path[fill=gray!35!white](-2.5,4.5)--++(1,0)--++(0,1)--++(-1,0);
\path[fill=gray!60!white](-3.5,5.5)--++(1,0)--++(0,-1)--++(1,0)--++(0,-1)--++(1,0)--++(0,-1)--++(3,0)--++(0,-1)--++(1,0)--++(0,-1)--++(1,0)--++(0,-1)--++(3,0)--++(0,-1)--++(-11,0);
\draw(-3,-1)node{$-30$};
\draw(-3,0)node{$-25$};
\draw(-3,1)node{$-20$};
\draw(-3,2)node{$-15$};
\draw(-3,3)node{$-10$};
\draw(-3,4)node{$-5$};
\draw(-3,5)node{$0$};
\draw(-2,-1)node{$-27$};
\draw(-2,0)node{$-22$};
\draw(-2,1)node{$-17$};
\draw(-2,2)node{$-12$};
\draw(-2,3)node{$-7$};
\draw(-2,4)node{$-2$};
\draw(-2,5)node{$3$};
\draw(-1,0)node{$-19$};
\draw(-1,1)node{$-14$};
\draw(-1,2)node{$-9$};
\draw(-1,3)node{$-4$};
\draw(-1,4)node{$1$};
\draw(-1,5)node{$6$};
%\draw(-1,6)node{$11$};
\draw(0,0)node{$-16$};
\draw(0,1)node{$-11$};
\draw(0,2)node{$-6$};
\draw(0,3)node{$-1$};
\draw(0,4)node{$4$};
\draw(0,5)node{$9$};
%\draw(0,6)node{$14$};
\draw(1,0)node{$-13$};
\draw(1,1)node{$-8$};
\draw(1,2)node{$-3$};
\draw(1,3)node{$2$};
\draw(1,4)node{$7$};
\draw(1,5)node{$12$};
%\draw(1,6)node{$17$};
\draw(2,0)node{$-10$};
\draw(2,1)node{$-5$};
\draw(2,2)node{$0$};
\draw(2,3)node{$5$};
\draw(2,4)node{$10$};
\draw(2,5)node{$15$};
%\draw(2,6)node{$20$};
\draw(3,0)node{$-7$};
\draw(3,1)node{$-2$};
\draw(3,2)node{$3$};
\draw(3,3)node{$8$};
\draw(3,4)node{$13$};
\draw(3,5)node{$18$};
%\draw(3,6)node{$23$};
\draw(4,0)node{$-4$};
\draw(4,1)node{$1$};
\draw(4,2)node{$6$};
\draw(4,3)node{$11$};
\draw(4,4)node{$16$};
\draw(4,5)node{$21$};
%\draw(4,6)node{$26$};
\draw(5,0)node{$-1$};
\draw(5,1)node{$4$};
\draw(5,2)node{$9$};
\draw(5,3)node{$14$};
\draw(5,4)node{$19$};
\draw(5,5)node{$24$};
%\draw(5,6)node{$29$};
\draw(6,0)node{$2$};
\draw(6,1)node{$7$};
\draw(6,2)node{$12$};
\draw(6,3)node{$17$};
\draw(6,4)node{$22$};
\draw(6,5)node{$27$};
%\draw(6,6)node{$32$};
\draw(7,0)node{$5$};
\draw(7,1)node{$10$};
\draw(7,2)node{$15$};
\draw(7,3)node{$20$};
\draw(7,4)node{$25$};
\draw(7,5)node{$30$};
%\draw(7,6)node{$35$};
\draw(-1,-1)node{$-24$};
\draw(0,-1)node{$-21$};
\draw(1,-1)node{$-18$};
\draw(2,-1)node{$-15$};
\draw(3,-1)node{$-12$};
\draw(4,-1)node{$-9$};
\draw(5,-1)node{$-6$};
\draw(6,-1)node{$-3$};
\draw(7,-1)node{$0$};
\draw[very thick](-1.5,4.5)--++(-1,0)--++(0,1)++(-1,0)--++(2,0)--++(0,-2)--++(2,0)--++(0,-1)++(-1,1)--++(0,-1)--++(4,0)--++(0,-1)++(-1,1)--++(0,-1)--++(1,0)--++(0,-1)--++(2,0)--++(0,-1)++(-1,1)--++(0,-1)--++(3,0)--++(0,-1);
\end{tikzpicture}
\]
We can encode this diagram by writing a B for each light-coloured box, with a sequence of Rs and Ds representing the path joining each box to the next. We see that we obtain a periodic sequence, with each period comprising $a$ Bs, $s-a$ Ds, and $t-a$ Rs. Conversely, any cyclic sequence of these symbols yields a bipartition $(\la,\mu)\in\bc sa\cap\bc ta$ in this way. For example, the diagram above corresponds to the cyclic sequence BDRBRR.

So we see that $\cardx{\bc sa\cap\bc ta}$ is the number of sequences comprising $a$ Bs, $s-a$ Ds, and $t-a$ Rs, modulo cyclic shifts. Counting these is a straightforward combinatorial exercise: disregarding cyclic shifts there are $(s+t-a)!/a!(s-a)!(t-a)!$ such sequences. None of these is fixed by any non-trivial cyclic shift, since the integers $a,s-a,t-a$ are coprime. So the final count is $(s+t-a-1)!/a!(s-a)!(t-a)!$.

\subsection{Extending Armstrong's Conjecture}

A recent exciting development in the theory of $(s,t)$-cores is Johnson's proof \cite{j} of Armstrong's Conjecture, which says that the average size of an $(s,t)$-core is $\frac1{24}(s-1)(t-1)(s+t+1)$. Of course, one can ask for the average size of a simultaneous core multipartition. Here we comment briefly on analogues of Armstrong's Conjecture for the two special cases mentioned in \cref{divsec,aasec}.

First consider the case where $s$ divides $t$; we enumerated simultaneous bicores in this situation in \cref{divsec}. We conjecture the average size of these bicores in two special subcases: where $s=t$, and where $t=0$.

\begin{conj}\label{ssarmstrong}
Suppose $0\ls a,b<s$, and that $s$ and $a-b$ are coprime. Then the average size of a bipartition in $\bc sa\cap\bc sb$ is
\[
\frac{(s+1)(a(s-a)+b(s-b)+1-s)}{12}.
\]
\end{conj}

\begin{eg}
\Yboxdim{7pt}
Take $s=5$, $a=1$ and $b=3$. Then $\cardx{\bc51\cap\bc53}=\frac15\binom51\binom53=10$. The Young diagrams of the ten bicores in $\bc51\cap\bc53$ are as follows.
\[
\begin{array}{c@{\qquad}c@{\qquad}c@{\qquad}c@{\qquad}c}
\varnothing\ \varnothing&\varnothing\ \yng(1)&\yng(1)\ \varnothing&\varnothing\ \yng(2)&\yng(1^2)\ \varnothing
\\[9pt]
\yng(1)\ \yng(1^2)&\yng(2)\ \yng(1)&\yng(1)\ \yng(3,1)&\yng(2,1^2)\ \yng(1)&\yng(3,1)\ \yng(2,1^2)
\end{array}
\]
These bipartitions have an average size of $3$, as predicted by \cref{ssarmstrong}.
\end{eg}

\begin{rmk}
In fact, Armstrong's Conjecture appears in disguise as a special case of \cref{ssarmstrong}. Suppose $a=0$. \ref{sposbetacore} shows that $(\la,\mu)\in\mc s{0,0}$ \iff $\la=\mu$ and $\la$ is an $s$-core. Applying \ref{sposbetacore} again, we find that $(\la,\la)\in\mc s{0,b}$ \iff $\la$ is both a $b$-core and an $(s-b)$-core. So in this case, $\smc\calt$ is simply the set of bipartitions $(\la,\la)$, where $\la$ is a $(b,s-b)$-core. (Note that such a partition is automatically an $s$-core.) So when $a=0$, \cref{ssarmstrong} is equivalent to Armstrong's Conjecture.
\end{rmk}

\begin{conj}\label{t0armstrong}
Suppose $0\ls a<s$ and $0\ls b$, and that $s$ and $a-b$ are coprime. Then the average size of a bipartition in $\bc sa\cap\bc 0b$ is
\[
\frac{(s+1)a(s-a)+(s-1)(b-1)(b+s+1)}{12}.
\]
\end{conj}

\begin{eg}
\Yboxdim{7pt}
Take $s=3$, $a=1$ and $b=2$. Then $\cardx{\bc31\cap\bc02}=\frac13\binom31\binom42=6$. The Young diagrams of the six bicores in $\bc31\cap\bc02$ are as follows.
\[
\varnothing\ \varnothing\qquad\varnothing\ \yng(1)\qquad\yng(1)\ \varnothing\qquad\varnothing\ \yng(2)\qquad\yng(1^2)\ \varnothing\qquad\yng(2)\ \yng(1^2)
\]
These bipartitions have an average size of $\frac53$, as predicted by \cref{t0armstrong}.
\end{eg}

Now consider the case $a=b$ addressed in \cref{aasec}. Here we make the following conjecture.

\begin{conj}\label{aaarmstrong}
Suppose $0\ls a<s\ls t$, and that $s$ and $t$ are coprime. Then the average size of a bipartition in $\bc sa\cap\bc ta$ equals
\[
\frac{(s-1)(t-1)(s+t-2a+1)-2a^2+2a}{12}.
\]
%Then there is a polynomial function $f_a(s,t)$ of degree $3$ (but of degree $2$ in each of $s,t$) with integer coefficients such that for $s,t$ coprime integers greater than $a$, the average size of a bipartition in $\bc sa\cap\bc ta$ equals $f_a(s,t)/12$.
%
%In particular, $f_0(s,t)=(s-1)(t-1)(s+t+1)$ and $f_1(s,t)=(s-1)(t-1)(s+t-1)$.
\end{conj}

Again, the case $a=0$ is equivalent to Armstrong's Conjecture, since $\bc s0\cap\bc t0$ is the set of bipartitions $(\la,\la)$ with $\la$ an $(s,t)$-core.

\begin{eg}
\Yboxdim{7pt}
Take $s=3$, $t=4$ and $a=1$. Then $\cardx{\bc31\cap\bc43}=\mfrac{5!}{1!2!3!}=10$. The Young diagrams of the ten bicores in $\bc31\cap\bc43$ are as follows.
\[
\begin{array}{c@{\qquad}c@{\qquad}c@{\qquad}c@{\qquad}c}
\varnothing\ \varnothing&\varnothing\ \yng(1)&\yng(1)\ \varnothing&\varnothing\ \yng(2)&\yng(1^2)\ \varnothing
\\[9pt]
\yng(1)\ \yng(1^2)&\yng(2)\ \yng(1)&\yng(2)\ \yng(1^2)&\yng(1^2)\ \yng(3,1^2)&\yng(3,1^2)\ \yng(2)
\end{array}
\]
These bipartitions have an average size of $3$, as predicted by \cref{aaarmstrong}.
\end{eg}

Johnson's proof of Armstrong's Conjecture relies on a geometric realisation of the set of $(s,t)$-cores, using Ehrhart theory. We hope to extend these ideas to core multipartitions in a future paper.

\end{document}